\documentclass[11pt]{article}
\usepackage{mathrsfs}
\usepackage{amsfonts}
\usepackage{bbm}
\usepackage{amsmath}
\usepackage{amssymb}
\usepackage[square, comma, sort&compress, numbers]{natbib}
\usepackage{xcolor}
\usepackage{tikz}
\usepackage{appendix}
\usepackage{lineno}
\usetikzlibrary{arrows,shapes,chains}
		
\newtheorem{theorem}{Theorem}[section]

\newtheorem{lemma}{Lemma}[section]

\newtheorem{remark}{Remark}[section]

\begin{document}

\title{Stochastic homogenization of nonlinear evolution equations with space-time nonlocality\thanks{The research is supported by the National Natural Science Foundation of China (No. 12171442). The research of J. Chen is partially supported by the CSC under grant No. 202206160033.}}
\author{Junlong Chen$^{1,}$\thanks{Email: chenjlm@hust.edu.cn (J. Chen)} and Yanbin Tang$^{1,2,}$\thanks{Corresponding author. Email: tangyb@hust.edu.cn (Y. Tang)}\\
\small {1) School of Mathematics and Statistics, Huazhong University of Science and}\\
\small {Technology, Wuhan, Hubei 430074, China}\\
\small {2) Hubei Key Laboratory of Engineering Modeling and Scientific Computing,}\\
\small {Huazhong University of Science and Technology, Wuhan, Hubei 430074, China}}

\date{\today}
\maketitle

{\bf Abstract} In this paper we consider the homogenization problem of nonlinear evolution equations with space-time non-locality, the problems are given by Beltritti and Rossi [JMAA, 2017, 455: 1470-1504]. When the integral kernel $J(x,t;y,s)$ is re-scaled in a suitable way and the oscillation coefficient $\nu(x,t;y,s)$ possesses periodic and stationary structure, we show that the solutions $u^{\varepsilon}(x,t)$ to the perturbed equations converge to $u_{0}(x,t)$, the solution of corresponding local nonlinear parabolic equation  as scale parameter $\varepsilon\rightarrow 0^{+}$. Then for the nonlocal linear index $p=2$ we give the convergence rate such that $||u^\varepsilon -u_{0}||_{_{L^{2}(\mathbb{R}^{d}\times(0,T))}}\leq C\varepsilon$. Furthermore, we obtain that the normalized difference $\frac{1}{\varepsilon}[u^{\varepsilon}(x,t)-u_{0}(x,t)]-\chi(\frac{x}{\varepsilon}, \frac{t}{\varepsilon^{2}}) \nabla_{x}u_{0}(x,t)$ converges to a solution of an SPDE with additive noise and constant coefficients. Finally, we give some numerical formats for solving non-local space-time homogenization.



\section{Introduction}\label{sec1}

Homogenization of local equations with rapidly oscillating coefficients both in time and in space has been extensively studied in recent years \cite{Holmbom1997,Pardoux2006,Piatnitski2015,Piatnitski2020,Shen2017,Shen2020,Armstrong,ChenTang2023}. The most of results are given for the equations whose coefficients are periodic and/or stationary. Homogenization of nonlinear parabolic equation on the qualitative results can be seen in \cite{Akagi2020,Akagi2021,Souganidis}. Geng and Shen \cite{Shen2017,Shen2020} considered periodic homogenization of parabolic equations with self-similar and non-self-similar scales in a bounded domain.

Piatnitski et al \cite{Piatnitski2015} studied the quantitative homogenization of equations
\begin{eqnarray*}
\partial_{t}u^{\varepsilon}(x,t)=\text{div}(a(\frac{x}{\varepsilon},\frac{t}{\varepsilon^{2}})\nabla u^{\varepsilon}(x,t))
\end{eqnarray*}
with space periodicity and time stationarity for the case of similar scales in the whole space, and they recently completed the results for the non-self-similar scales \cite{Piatnitski2020}. For the homogenization of nonlocal operators, we refer to \cite{Piatnitski2017,PiatnitskiStochastic} and references therein.

In this paper we consider the homogenization problem of nonlinear evolution equations with space-time non-locality, which are given by Beltritti and Rossi \cite{Rossi2017} for $p\geq2$ and Aimar, Beltritti and Gomez \cite{Aimar2018} for $p=2$,
\begin{eqnarray}\label{Equ1.1}
\left\{\begin{array}{l}
\mathbb{S}(J,\nu,v,g)\equiv\displaystyle{\iint_{\mathbb{R}^{d}\times\mathbb{R}}J(x-y,t-s)\nu(x,t;y,s)a(\bar{v}(y,s)-v(x,t))dyds=0,}\\
a(\bar{v}(y,s)-v(x,t))=|\bar{v}(y,s)-v(x,t)|^{p-2}(\bar{v}(y,s)-v(x,t)),\end{array}\right.
\end{eqnarray}
for $(x,t)\in\mathbb{R}^{d}\times[0,\infty),$ $p\geq 2,$ $\nu(x,t;y,s)$ is a bounded measurable function, $\bar{v}(x,t)$ stands for an extension of a given function $g(x,t)$, that is,
\begin{eqnarray}\label{Equ1.2}
\bar{v}(x,t)=\left\{\begin{array}{ll}v(x,t), &t>0,\\ g(x,t), &t\leq0. \end{array}\right.
\end{eqnarray}
The kernel $J(x,t): \mathbb{R}^{d+1}\rightarrow\mathbb{R}$ is radially symmetric in $x$ and satisfies that
\begin{eqnarray}\label{Equ1.3}
J\geq 0,\; \iint_{\mathbb{R}^{d}\times\mathbb{R}}J(x,t)dxdt=1.
\end{eqnarray}
We focus on the homogenization of the nonlocal equation (\ref{Equ1.1}) and hope to obtain qualitative homogenization of equations with mixed space-time periodic and/or stationary structure and then give the convergence rate of the solutions.

The fractional order nonlocal parabolic operators have recently attracted great interest in partial differential equation (PDE) and homogenization theory. For example, Stinga and Torrea \cite{Torrea} considered the regularity of solution to space-time nonlocal equations driven by the fractional power of the heat operator $H=\partial_t-\Delta$, i.e., $H^{s}=(\partial_t-\Delta)^s (0<s<1)$, it has an integral form
\begin{eqnarray}
&&H^{s}u(x,t)=\int_{0}^{\infty}\int_{\mathbb{R}^{d}}N_{s}(y,\tau)\Big(u(x,t)-u(x-y,t-\tau)\Big)dyd\tau, \label{Equ1.4}\\
&&N_{s}(y,\tau)=\frac{s}{(4\pi)^{d/2}\Gamma(1-s)\tau^{d/2+1+s}}e^{-\frac{|y|^{2}}{4\tau}}, \label{Equ1.5}
\end{eqnarray}
where $N_{s}(x,t)$ represents the nonlocal (non-integrable) kernel, it differs from the integrable kernel $J$ in (\ref{Equ1.1}) and (\ref{Equ1.3}). Some homogenization results may arise when the kernel $N_{s}(x,t)$ is re-scaled in a suitable way, we will discuss this structure in subsequent research.

To the best of our knowledge, no one has investigated quantitative homogenization of nonlinear nonlocal parabolic equation with space-time coefficients, even in the period environment. In this paper we first consider nonlinear nonlocal parabolic equation (\ref{Equ1.1}) with space-time mixed coefficients. Furthermore one can also consider the general form of equation (\ref{Equ1.1}) such that $a(\cdot)$ is strongly monotone and Lipschitz continuous. Subsequent research will be concentrated the structure $\nu(x,t;\frac{x}{\varepsilon},\frac{t}{\varepsilon^r})(0<r<+\infty)$ under the non-self-similar scale condition.

Let $g:\mathbb{R}^{d+1}\rightarrow\mathbb{R}$ be a bounded function such that $g$ is $C^{2}$ with bounded derivatives. Denote $\overline{\mathcal{C}}$ the set of uniform continuous functions and $L^{\infty}(g)$ stands for the set of bounded functions with norm less or equal than $||g||_{L^\infty (R^d\times(-\infty,0))}$. Assume that
\begin{eqnarray}\label{Equ1.6}
\nu(x,t;y,s)=\lambda(x,t)\mu(y,s),\; 0<\alpha_1\leq \lambda,\mu\leq\alpha_2.
\end{eqnarray}
To highlight our ideas, we choose the simplest structure $\nu(x,t;y,s)=\mu(x,t)\mu(y,s)$ and actually $\nu$ can be extended to positive function (no symmetry required).

For positive constants $\delta$ and $\gamma$ we denote the set
\begin{eqnarray}\label{Equ1.7}
D_{\delta\gamma}=\{(x,t)\in\mathbb{R}^{d+1}:\delta\leq t\leq \delta+\gamma\}.
\end{eqnarray}
The result of Beltritti and Rossi \cite{Rossi2017} reads as follows.
\begin{theorem}\cite{Rossi2017}\label{The1.1} ({\bf Uniqueness and existence}) Let the kernel $J(x,t)$ be nonnegative, continuous and compactly supported in the set $\{(x,t)\in\mathbb{R}^{d+1}:t\geq 0\}$ with $\iint_{\mathbb{R}^{d}\times\mathbb{R}}J(x,t)dxdt=1$. If $g\in L^{\infty}(\mathbb{R}^{d}\times(-\infty,0))$, there exists a unique function $u\in\overline{\mathcal{C}}\cap L^{\infty}(g)(\mathbb{R}^{d}\times[0,+\infty))$ such that $u$ solves the equation (\ref{Equ1.1}): $\mathbb{S}(J,1,u,g)=0$.\\
({\bf Scaling invariance}) Let $\varepsilon>0$ and $J_{\varepsilon}(x,t)=\frac{1}{\varepsilon^{d+2}}J(\frac{x}{\varepsilon},\frac{t}{\varepsilon^{2}})$. If $u^{\varepsilon}(x,t)$ is a solution to the equation $\mathbb{S}(J_{\varepsilon},1,u^{\varepsilon},g)=0$ and $v$ solves $\mathbb{S}(J,1,v,\varepsilon^{d+2}g(\varepsilon x,\varepsilon^{2}t))=0$, then $u^{\varepsilon}(x,t)=\frac{1}{\varepsilon^{d+2}}v(\frac{x}{\varepsilon},\frac{t}{\varepsilon^{2}}).$\\
({\bf Convergence}) If $J(x,t)$ satisfies (\ref{Equ1.3}) and is compactly supported in $D_{\delta\gamma}$ and $f:\mathbb{R}^{d}\rightarrow \mathbb{R}$ is a bounded function such that it is $C^{2}$ with bounded derivatives. As $\varepsilon\rightarrow 0^{+}$, the solutions $u^{\varepsilon}(x,t)$ to $\mathbb{S}(J_{\varepsilon},1,u^{\varepsilon},f)=0$ converge up to a subsequence to a viscosity solution to the nonlinear degenerate parabolic PDE: $C(J)||\nabla u||^{p-2}\frac{\partial u}{\partial t}=\Delta_{p}u$ with initial condition $u(x,0)=f(x)$, where $\Delta_{p}$ is $p-$Laplacian operator and $C(J)$ is a constant depending on $J$.
\end{theorem}

When the integral kernel $J(x,t)$ is re-scaled in a suitable way and the oscillation coefficient $\nu(x,t)$ possesses periodic and stationary structure, we show that the solutions $u^{\varepsilon}(x,t)$ to the perturbed equations converge to $u_{0}(x,t)$, the solution of the corresponding local nonlinear partial differential equations as $\varepsilon\rightarrow 0^{+}$. Then for the nonlocal linear index $p=2$ we give the convergence rate such that $||u^\varepsilon-u_{0}||_{_{L^{2}(\mathbb{R}^{d}\times(0,T))}}\leq C\varepsilon$. Furthermore, we obtain that the normalized difference $\frac{1}{\varepsilon}[u^{\varepsilon}(x,t)-u_{0}(x,t)]-\chi(\frac{x}{\varepsilon},\frac{t}{\varepsilon^{2}})\nabla_{x}u_{0}$ converges to a solution of a stochastic partial differential equation (SPDE) with additive noise and constant coefficients.

An interesting conclusion  in \cite[Theorem 3.1]{Harlim} compared to our article when $\nu(x,y)=\mu(x)\mu(y)$ does not depend on $\varepsilon$, and have the following estimate $||u^\varepsilon-u||_{L^{2}}\leq C\varepsilon^\epsilon$, where $\epsilon$ can be arbitrarily close to $1$, but less than $1$. In fact, this is due to Taylor expand we can collect all $\varepsilon^1-$terms, first-order moment of kernel function which is odd function and equal to zero. But integral of $\varepsilon^1-$term is no longer equal to zero when $\nu(x,y)=\mu(\frac{x}{\varepsilon})\mu(\frac{y}{\varepsilon}).$

\section{Statement of problem and main results}\label{sec2}
\setcounter{equation}{0}

We consider the following nonlocal scaling $p-$Laplacian ($p\geq 2$) operator with space-time kernel
\begin{eqnarray}\label{Equ2.1}
\mathbb{S}(J_{\varepsilon},\nu,u^{\varepsilon},g)&\equiv&\frac{1}{\varepsilon^{d+p+2}}\int\int_{\mathbb{R}^{d}\times\mathbb{R}}J(\frac{x-y}{\varepsilon},\frac{t-s}{\varepsilon^2}) \nu(\frac{x}{\varepsilon},\frac{t}{\varepsilon^2};\frac{y}{\varepsilon},\frac{s}{\varepsilon^2})\nonumber\\
&&{\hspace{-1cm}}\cdot\Big|\bar{u}^{\varepsilon}(y,s)-u^{\varepsilon}(x,t)\Big|^{p-2}\Big(\bar{u}^{\varepsilon}(y,s)-u^{\varepsilon}(x,t)\Big)dyds=0,
\end{eqnarray}
where the space-time coefficient $\nu(x,t;y,s)=\mu(x,t)\mu(y,s)$ satisfies (\ref{Equ1.6}).

In this paper we will consider three cases of $\mu(x,t)$ as follows.
\begin{eqnarray}\label{Equ2.2}
\left\{\begin{array}{l}
{\bf Case~ 1}: \mu(x,t) ~is~ \mathbb{Z}^d\times\mathbb{Z} -periodic ~ in ~x,t , \\
{\bf Case~ 2}: \mu(x,t) ~is~ \mathbb{Z}^d -periodic ~ in ~x ~ and ~ \mathbb{R} -stationary ~ in~t,\\
{\bf Case~ 3}: \mu(x,t) ~is~ \mathbb{R}^d\times\mathbb{R} -stationary ~ in ~x,t. \end{array}\right.
\end{eqnarray}

The stationarity of the hypothesis in this paper means that  $\nu=\mu(x,t)\mu(y,s)$ with a statistically homogeneous ergodic field $\mu$ which assembles with  an ergodic group of measurable transformations.  Specifically, $T$ is a  (space or time) dynamical system, $(\Omega,\digamma,P)$ be a standard probability space and we assume that $\mu(x,w)=\mu(T_{x}w),T_x:\Omega\to\Omega$ satisfy the following properties.

\noindent(1)\quad $T_{y_1}T_{y_2}=T_{y_1+y_2}$ for all $y_1,y_2$ in $\mathbb{R}^d,T_0=Id$.

\noindent(2)\quad $P(T_{y}A)=P(A)$ for all $A\in\mathcal{F}$ and all $y\in \mathbb{R}^d$.

\noindent(3)\quad $T_x$ is measurable map from $\mathbb{R}^d\times\Omega$ to $\Omega$, where $\mathbb{R}^d$ is equipped with the Borel $\sigma-$algebra.

We focus on the {\bf Case 1} and {\bf Case 2}. Specially for the linear case $p=2$ in (\ref{Equ2.1}) we get the normalized difference between $u^\varepsilon(x,t)$ and the first two terms of asymptotic expansion and prove that it converges to a solution of SPDE in {\bf Case 2}. For the {\bf Case 3} we just give the results and omit the proof in details. We now state our main results.
\subsection{$\mu(x,t)$ is periodic in $x,t$}\label{sec21}
\begin{theorem}\label{The2.1} $J(x,t):\mathbb{R}^{d+1}\rightarrow\mathbb{R}$ is compactly supported in $D_{\delta\gamma}$ defined in (\ref{Equ1.7}) and $f:\mathbb{R}^{d}\rightarrow \mathbb{R}$ is a bounded function such that it is $C^{2}$ with bounded derivatives. The solutions $u^{\varepsilon}(x,t)$ to $\mathbb{S}(J_{\varepsilon},\nu,u^{\varepsilon},f)=0$ converge along subsequence uniformly as $\varepsilon \to 0$ on compact sets to a viscosity solution $u(x,t)$  as $\varepsilon\rightarrow 0^{+}$, such that
\begin{eqnarray}\label{Equ2.3}
\left\{\begin{array}{l}
\mathfrak{N}(\nabla u)\frac{\partial u}{\partial t}=\mathfrak{P}(\nabla u,\nabla\nabla u),\\ u(x,0)=f(x), \end{array}\right.
\end{eqnarray}
where $\mathfrak{N}$ and $\mathfrak{P}$ are given by
\begin{eqnarray}\label{Equ2.4}
&&\mathfrak{N}(\nabla u)=\int_{\mathbb{T}^d}\int_{\mathbb{T}}\int_{\mathbb{R}}\int_{\mathbb{R}^d}J(z,r)\mu(\xi,q)\mu(\xi-z,q-r)\nonumber\\
&&\quad\cdot\Big|\nabla u(x,t)\Big(-z+\chi_1(\xi-z,q-r)-\chi_1(\xi,q)\Big)\Big|^{p-2}rdzdrd\xi dq,
\end{eqnarray}
\begin{eqnarray}\label{Equ2.5}
&&\mathfrak{P}(\nabla u,\nabla\nabla u)=\int_{\mathbb{T}^d}\int_{\mathbb{T}}\int_{\mathbb{R}}\int_{\mathbb{R}^d}J(z,r)\mu(\xi,q)\mu(\xi-z,q-r)\nonumber\\
&&\qquad\cdot\Big|\nabla u(x,t)\Big(-z+\chi_1(\xi-z,q-r)-\chi_1(\xi,q)\Big)\Big|^{p-2}\nonumber\\
&&\qquad\cdot\Big[-z\otimes\chi_1(\xi-z,q-r)+{\frac{1}{2}}z\otimes z\Big]\nabla\nabla u(x,t)dzdrd\xi dq.
\end{eqnarray}
\end{theorem}

Due to the classical method of asymptotic expansion \cite{BensoussanLions}, we first construct some auxiliary functions to prove our main results in Theorem \ref{The2.1}.

\subsection{Correctors and auxiliary cell problems}\label{sec22}
\begin{lemma}\label{lem2.1} If there exist two periodic functions $\chi_1\in C_b(\mathbb{T}^d\times\mathbb{T})^d, \chi_2\in C_b(\mathbb{T}^d\times\mathbb{T})^{d^2}$  such that
\begin{eqnarray}\label{Equ2.6}
w^\varepsilon(x,t)=u(x,t)+\varepsilon\chi_1(\frac{x}{\varepsilon},\frac{t}{\varepsilon^2})\nabla u(x,t)+\varepsilon^2\chi_2(\frac{x}{\varepsilon},\frac{t}{\varepsilon^2})\nabla \nabla u(x,t),
\end{eqnarray}
then as $\varepsilon \to 0^{+}$ we have
\begin{eqnarray}\label{Equ2.7}
\mathbb{S}(J_{\varepsilon},\nu,w^{\varepsilon},f)=\frac{1}{\varepsilon}r_0(x,t)+r_{\varepsilon}(x,t)+\Phi^*_{\varepsilon}(x,t)+(p-2)o(\varepsilon),
\end{eqnarray}
where
\begin{eqnarray}\label{Equ2.8}
r_0(u)(x,t)&=&r_0(x,t)=\int_{\mathbb{R}}\mu(\frac{x}{\varepsilon},\frac{t}{\varepsilon^2})\int_{\mathbb{R}^d}J(z,r)\mu(\frac{x}{\varepsilon}-z,\frac{t}{\varepsilon^2}-r)\nonumber\\
&&{\hspace{-1cm}}\cdot\Big|\Big(-z+\chi_1(\frac{x}{\varepsilon}-z,\frac{t}{\varepsilon^2}-r)-\chi_1(\frac{x}{\varepsilon},\frac{t}{\varepsilon^2})\Big)\cdot\nabla u(x,t)\Big|^{p-2}\nonumber\\
&&{\hspace{-1cm}}\cdot\Big[-z+\chi_1(\frac{x}{\varepsilon}-z,\frac{t}{\varepsilon^2}-r)-\chi_1(\frac{x}{\varepsilon},\frac{t}{\varepsilon^2})\Big]\cdot\nabla u(x,t)dzdr,
\end{eqnarray}
\begin{eqnarray}\label{Equ2.9}
r_{\varepsilon}(u)(x,t)&=&r_{\varepsilon}(x,t)=(p-1)\int_{\mathbb{R}}\mu(\frac{x}{\varepsilon},\frac{t}{\varepsilon^2})\int_{\mathbb{R}^d}J(z,r)\mu(\frac{x}{\varepsilon}-z,\frac{t}{\varepsilon^2}-r)\nonumber\\
&&{\hspace{-1cm}}\cdot\Big|\Big(-z+\chi_1(\frac{x}{\varepsilon}-z,\frac{t}{\varepsilon^2}-r)-\chi_1(\frac{x}{\varepsilon},\frac{t}{\varepsilon^2})\Big)\cdot\nabla u(x,t)\Big|^{p-2}\nonumber\\
&&{\hspace{-1cm}}\cdot\Big\{\Big[-z\otimes \chi_1(\frac{x}{\varepsilon}-z,\frac{t}{\varepsilon^2}-r)+\frac{1}{2}z^{2}\Big]\cdot\nabla\nabla u\nonumber\\
&&{\hspace{-1cm}}+\Big[\chi_2(\frac{x}{\varepsilon}-z,\frac{t}{\varepsilon^2}-r)-\chi_2(\frac{x}{\varepsilon},\frac{t}{\varepsilon^2})\Big]\cdot\nabla\nabla u-r\frac{\partial u(x,t)}{\partial t}\Big\}dzdr,
\end{eqnarray}
\begin{eqnarray}\label{Equ2.10}
&&\Phi^*(u)(x,t)=\Phi^*_{\varepsilon}(x,t)=(p-1)\int_{\mathbb{R}}\mu(\frac{x}{\varepsilon},\frac{t}{\varepsilon^2})\int_{\mathbb{R}^{d}}J(z,r)\mu(\frac{x}{\varepsilon}-z,\frac{t}{\varepsilon^2}-r) \nonumber\\
&&\qquad\qquad\qquad\cdot\Big|\Big(-z+\chi_1(\frac{x}{\varepsilon}-z,\frac{t}{\varepsilon^2}-r)-\chi_1(\frac{x}{\varepsilon},\frac{t}{\varepsilon^2})\Big)\cdot \nabla u(x,t)\Big|^{p-2}\nonumber\\
&&\qquad\qquad\qquad\cdot\Big\{\Big[\int_{0}^{1}\nabla\nabla u(x-\varepsilon z\theta,t-\varepsilon^2 r\theta)\cdot z\otimes z(1-\theta)d\theta\nonumber\\
&&-\frac{1}{2}\nabla\nabla u(x,t)\cdot z\otimes z+2\varepsilon r\int_{0}^{1}\partial_t\nabla u(x-\varepsilon z\theta,t-\varepsilon^2 r\theta)\otimes z(1-\theta) d \theta\nonumber\\
&&+\varepsilon\chi_{1}\Big(\frac{x}{\varepsilon}-z,\frac{t}{\varepsilon^2}-r\Big)\int_{0}^{1}\nabla\nabla\nabla u(x-\varepsilon z\theta,t-\varepsilon^2 r\theta)z\otimes z(1-\theta)d\theta \nonumber\\
&&-\varepsilon\chi_{2}\Big(\frac{x}{\varepsilon}-z,\frac{t}{\varepsilon^2}-r\Big)\cdot\int_{0}^{1}\nabla\nabla\nabla u(x-\varepsilon z\theta,t-\varepsilon^2 r\theta)zd\theta\Big]\nonumber\\
&&-\varepsilon\Big[r\chi_1(\frac{x}{\varepsilon}-z,\frac{t}{\varepsilon^2}-r)\cdot\frac{\partial\nabla u}{\partial t}\Big]
+\varepsilon^2 r^2\int_{0}^{1}\partial_{tt} u(x-\varepsilon z\theta,t-\varepsilon^2 r\theta)(1-\theta) d \theta\nonumber\\
&&+2\varepsilon^2\chi_1\Big(\frac{x}{\varepsilon}-z,\frac{t}{\varepsilon^2}-r\Big)\int_{0}^{1}\partial_{t}\nabla\nabla  u(x-\varepsilon z\theta,t-\varepsilon^2r\theta)\otimes z r (1-\theta)d\theta\nonumber\\
&&-\varepsilon^2\chi_{2}\Big(\frac{x}{\varepsilon}-z,\frac{t}{\varepsilon^2}-r\Big)\int_{0}^{1}\partial_t\nabla\nabla u(x-\varepsilon z\theta,t-\varepsilon^2 r\theta)rd\theta\nonumber\\
&&+\varepsilon^3\chi_1\Big(\frac{x}{\varepsilon}-z,\frac{t}{\varepsilon^2}-r\Big)\int_{0}^{1}\partial_{tt}\nabla u(x-\varepsilon z\theta,t-\varepsilon^2r\theta)r^2 (1-\theta)d\theta\Big\}dzdr,
\end{eqnarray}
where $z^{2}=z\otimes z=\{z^{i}z^{j}\}^{d}_{i,j=1}$ and
\begin{eqnarray*}
\nabla\nabla u(\cdot)z=\frac{\partial^{2}u}{\partial x^{i}\partial x^{j}}(\cdot)z^{j},\;\nabla\nabla\nabla u(\cdot)z\otimes z=\frac{\partial^{3}u}{\partial x^{i}\partial x^{j}\partial x^{k}}(\cdot)z^{j}z^{k}.
\end{eqnarray*}
\end{lemma}
\noindent{\bf Proof.} Similar to the proof of $\|\phi_{\varepsilon}\|_{L^{2}(\mathbb{R}^{d})}\rightarrow 0$ as $\varepsilon \rightarrow 0^{+}$ in \cite[Proposition 5]{Piatnitski2017} we can get $\Phi^*_{\varepsilon}(x,t)\rightarrow 0$ as $\varepsilon \rightarrow 0^{+}$. Set $\frac{x-y}{\varepsilon}=z,\frac{t-s}{\varepsilon^2}=r$, we get
\begin{eqnarray}\label{Equ2.11}
&&\mathbb{S}(J_{\varepsilon},\nu,w^{\varepsilon},f)\nonumber\\
&=&\int_{\mathbb{R}}\int_{\mathbb{R}^d}\frac{1}{\varepsilon^{d+2+p}}J(z,r)\mu(\frac{x}{\varepsilon},\frac{t}{\varepsilon^2})\mu(\frac{x}{\varepsilon}-z,\frac{t}{\varepsilon^2}-r)
\Big|u(x-\varepsilon z,t-\varepsilon^2 r)\nonumber\\
&+&\varepsilon\chi_1(\frac{x}{\varepsilon}-z,\frac{t}{\varepsilon^2}-r)\nabla u(x-\varepsilon z)+\varepsilon^2\chi_2(\frac{x}{\varepsilon}-z,\frac{t}{\varepsilon^2}-r)\nabla\nabla u(x-\varepsilon z,t-\varepsilon^2 r)\nonumber\\
&-&u(x,t)-\varepsilon\chi_1(\frac{x}{\varepsilon},\frac{t}{\varepsilon^2})\nabla u(x,t)-\varepsilon^2\chi_2(\frac{x}{\varepsilon},\frac{t}{\varepsilon^2})\nabla  \nabla u(x,t)\Big|^{p-2}\Big[u(x-\varepsilon z,t-\varepsilon^2 r)\nonumber \\
&+&\varepsilon\chi_1(\frac{x}{\varepsilon}-z,\frac{t}{\varepsilon^2}-r)\nabla u(x-\varepsilon z,t-\varepsilon^2 r)-u(x,t)-\varepsilon\chi_1(\frac{x}{\varepsilon},\frac{t}{\varepsilon^2})\nabla u(x,t)\nonumber\\
&+&\varepsilon^2\chi_2(\frac{x}{\varepsilon}-z,\frac{t}{\varepsilon^2}-r)\nabla \nabla u(x-\varepsilon z,t-\varepsilon^2 r)-\varepsilon^2\chi_2(\frac{x}{\varepsilon},\frac{t}{\varepsilon^2})\nabla\nabla u(x,t)\Big]dzdr.
\end{eqnarray}
Using the Taylor expansion formula with integral remainder
\begin{eqnarray}\label{Equ2.12}
\left\{\begin{array}{l}
u(y)=u(x)+\displaystyle{\int_{0}^{1}\frac{\partial}{\partial \theta}u(x+(y-x)\theta)d\theta}=u(x)+\displaystyle{\int_{0}^{1} \nabla u(x+(y-x)\theta)\cdot(y-x)d\theta},\\
u(y)=u(x)+\nabla u(x)\cdot(y-x)+\displaystyle{\int_{0}^{1}\nabla\nabla u(x+(y-x)\theta)(y-x)\cdot(y-x)(1-\theta)d\theta}, \end{array}\right.
\end{eqnarray}
we have
\begin{eqnarray}\label{Equ2.13}
&&u(x-\varepsilon z,t-\varepsilon^2 r)-u(x,t)+\varepsilon \nabla u(x,t)\cdot z+\varepsilon^{2}\frac{\partial u(x,t)}{\partial t}r\nonumber\\
&&=\frac{\varepsilon^{2}}{2}\sum_{i,j=1}^{n}\frac{\partial^{2}u(x,t)}{\partial x_{i}\partial x_{j}}z_{i}z_{j}+o(\varepsilon^{2}).
\end{eqnarray}
Putting (\ref{Equ2.12}) and (\ref{Equ2.13}) into (\ref{Equ2.11}), we have (\ref{Equ2.7}), where the functions $r_0(x,t),$ $r_{\varepsilon}(x,t),$ $\Phi^*_{\varepsilon}(x,t)$ are defined in (\ref{Equ2.8}), (\ref{Equ2.9}) and (\ref{Equ2.10}) respectively. Due to the order of $\varepsilon$, we put the remainder terms into the fourth part, that is, it contains many higher-order terms in $\varepsilon$. For the given functions $\chi_1$ and $\chi_2$, we need to show that (1) $r_0(x,t)=0$, (2) $r_{\varepsilon}(x,t)$ converges in the sense of viscous solution  as $\varepsilon \rightarrow 0^{+}$ to a solution of a local nonlinear equation, and (3) $\Phi^*_{\varepsilon}(x,t)\to 0$, (4) the fourth part is an infinitesimal $o(\varepsilon)$ as $\varepsilon \to 0^{+}$. This completes the proof of Lemma \ref{lem2.1}. $\Box$

The formula (\ref{Equ2.7}) gives an asymptotic decomposition of $\mathbb{S}(J_{\varepsilon},\nu,w^{\varepsilon},f)$ in $\varepsilon$, we now deal with the first three parts $r_0(x,t), r_{\varepsilon}(x,t), \Phi^*_{\varepsilon}(x,t)$ and get more precisely asymptotic behavior.
\begin{lemma}\label{lem2.2} For the function $u(x,t)$ given in (\ref{Equ2.6}), the equation $r_0(x,t)=0$ has a periodic 
	solution  $\chi_1(\xi,q)\in  C_b(\mathbb{T}^d\times \mathbb{T}),q>0; \chi_1(\xi,q)=0, q\leq0$ .
\end{lemma}
\noindent{\bf Proof.} Denote $\hbar(\xi,q)=\chi_1(\xi,q)+\xi-\frac{1}{2},$ we have
\begin{eqnarray}\label{Equ2.14}
	&&\int_{\mathbb{R}}\int_{\mathbb{R}^d}J(z,r)\mu(\xi-z,q-r)\Big|\nabla u(x,t)\cdot\Big(-z+\chi_1(\xi-z,q-r)-\chi_1(\xi,q)\Big)\Big|^{p-2}\nonumber\\
	&&\quad\cdot\Big(\nabla u(x,t)\cdot\Big(-z+\chi_1(\xi-z,q-r)-\chi_1(\xi,q)\Big)\Big)dzdr\nonumber\\
	&&=\int_{\mathbb{R}}\int_{\mathbb{R}^d}J(z,r)\mu(\xi-z,q-r)\Big|\nabla u(x,t)\cdot\Big(\hbar(\xi-z,q-r)-\hbar(\xi,q)\Big)\Big|^{p-2}\nonumber\\
	&&\quad\cdot\Big(\nabla u(x,t)\cdot\Big(\hbar(\xi-z,q-r)-\hbar(\xi,q)\Big)\Big)dzdr:=C\hbar.
\end{eqnarray}
We can obtain the existence of function $\hbar(\xi,q) \in \overline{\mathcal{C}} $ from \cite[Theorem1.1]{Rossi2017}  . Furthermore there is a function $\chi_1(\xi,q)=\chi_1(\frac{x}{\varepsilon}, \frac{t}{\varepsilon^2})\in C_b(\mathbb{T}^d\times \mathbb{T})$ such that $r_0(x,t)=0$.
$\Box$

For the function $\chi_1(\xi,q)=\chi_1(\frac{x}{\varepsilon},\frac{t}{\varepsilon^2})$ determined in Lemma \ref{lem2.1}, we denote
\begin{eqnarray}
&&\mathfrak{J}_{0}(z,r,x,t,\frac{x}{\varepsilon},\frac{t}{\varepsilon^2})=J(z,r)\mu(\frac{x}{\varepsilon},\frac{t}{\varepsilon^2})\mu(\frac{x}{\varepsilon}-z,\frac{t}{\varepsilon^2}-r)\nonumber\\
&&\qquad\qquad\cdot\Big|\Big(-z+\chi_1(\frac{x}{\varepsilon}-z,\frac{t}{\varepsilon^2}-r)-\chi_1(\frac{x}{\varepsilon},\frac{t}{\varepsilon^2})\Big)\nabla u(x,t)\Big|^{p-2},\label{Equ2.15}\\
&&\Re_{0}(x,t,\frac{x}{\varepsilon},\frac{t}{\varepsilon^2})=\int_{\mathbb{R}}\int_{\mathbb{R}^d}\mathfrak{J}_{0}(z,r,x,t,\frac{x}{\varepsilon},\frac{t}{\varepsilon^2})\nonumber\\
&&\qquad\qquad\cdot\Big[-z\otimes\chi_1(\frac{x}{\varepsilon}-z,\frac{t}{\varepsilon^2}-r)+\frac{1}{2}z\otimes z\Big]\nabla\nabla u(x,t)dzdr,\label{Equ2.16}
\end{eqnarray}
and introduce an operator $\mathcal{A}$ such that
\begin{eqnarray}\label{Equ2.17}
&&(\mathcal{A}\chi_2)(x,t,\frac{x}{\varepsilon},\frac{t}{\varepsilon^2}):=\int_{\mathbb{R}}\int_{\mathbb{R}^d}\mathfrak{J}_{0}(z,r,x,t,\frac{x}{\varepsilon},\frac{t}{\varepsilon^2})\nonumber\\
&&\qquad\qquad\cdot\Big[\chi_2(\frac{x}{\varepsilon}-z,\frac{t}{\varepsilon^2}-r)-\chi_2(\frac{x}{\varepsilon},\frac{t}{\varepsilon^2})\Big]\nabla\nabla udzdr.
\end{eqnarray}
\begin{lemma}\cite{Piatnitski2017}\label{lem2.3} The linear equation $(\mathcal{A}\chi_2)(x,t,\xi,q)=\mathfrak{P}(\nabla u,\nabla\nabla u)-\Re_{0}(x,t,\xi,q)$ has a periodic solution $\chi_2(\xi,q)\in  C_b(\mathbb{T}^d\times \mathbb{T}),q>0; \chi_2(\xi,q)=0, q\leq0$ .
\end{lemma}

Using symmetry and periodicity and \cite[Proposition 7]{Piatnitski2017} we have
\begin{eqnarray}\label{Equ2.18}
\int_{\mathbb{T}}\int_{\mathbb{T}^d}(\mathcal{A}\chi_2)(x,t,\xi,q)d\xi dq=0.
\end{eqnarray}
Following the previous Lemmas \ref{lem2.1}-\ref{lem2.3} we can give the proof of Theorem \ref{The2.1}.

\subsection{Proof of Theorem \ref{The2.1}}\label{sec23}

\noindent{\bf Proof of Theorem \ref{The2.1}.} By Taylor's theorem, for $|b|<|a|$ we have
\begin{eqnarray}\label{Equ2.19}
|a+b|^{p-2}(a+b)=|a|^{p-2} a+(p-1)|a|^{p-2} b+(p-1)(p-2)|\eta|^{p-4}\eta b^{2},
\end{eqnarray}
where $\eta$ lies in the segment that joins $a$ with $b$, and
\begin{eqnarray}\label{Equ2.20}
&&u(x-\varepsilon z,t-\varepsilon^{2}r)-u(x,t)=-\varepsilon\nabla u(x,t)\cdot z\nonumber \\
&&\qquad-\varepsilon^{2}r\partial_{t}u(x,t)+\frac{1}{2}\varepsilon^{2}\sum_{i,j=1}^{d}\partial_{x_{i}x_{j}}u(x,t)z_{i}z_{j}+o(\varepsilon^{2}).
\end{eqnarray}
Denote $\xi=\frac{x}{\varepsilon}, q=\frac{t}{\varepsilon^2}$ on the period $\xi\in \mathbb{T}^d, q\in \mathbb{T},$ using (\ref{Equ2.19}) and (\ref{Equ2.20}), we have
\begin{eqnarray}\label{Equ2.21}
r_0(x,t)&=&\int_{\mathbb{R}}\mu(\xi,q)\int_{\mathbb{R}^d}J(z,r)\mu(\xi-z,q-r)\nonumber\\
&&\cdot\Big|\Big(-z+\chi_1(\xi-z,q-r)-\chi_1(\xi,q)\Big)\nabla u(x,t)\Big|^{p-2}\nonumber\\
&&\cdot\Big[-z+\chi_1(\xi-z,q-r)-\chi_1(\xi,q)\Big]\nabla u(x,t)dzdr,
\end{eqnarray}
\begin{eqnarray}\label{Equ2.22}
r_{\varepsilon}(x,t)&=&(p-1)\int_{\mathbb{R}}\mu(\xi,q)\int_{\mathbb{R}^d}J(z,r)\mu(\xi-z,q-r)\nonumber\\
&&\cdot\Big|\Big(-z+\chi_1(\xi-z,q-r)-\chi_1(\xi,q)\Big)\nabla u(x,t)\Big|^{p-2}\nonumber\\
&&\cdot\Big\{\Big[-z\otimes \chi_1(\xi-z,q-r)+{\frac{1}{2}}z\otimes z\Big]\nabla\nabla u\nonumber\\
&&+\Big[\chi_2(\xi-z,q-r)-\chi_2(\xi,q)\Big]\nabla\nabla u-r\frac{\partial u(x,t)}{\partial t}\Big\}dzdr:=I^{\varepsilon}_{1}+I^{\varepsilon}_2.
\end{eqnarray}
Lemma \ref{lem2.2} implies that $r_0(x,t)=0$, and from Lemma \ref{lem2.3}, $\mu$ and $\chi_1, \chi_2$ will operate on the torus. 

We have
\begin{eqnarray}\label{Equ2.23}
\lim_{\varepsilon\rightarrow 0}I_1^{\varepsilon}&=&(p-1)\int_{\mathbb{T}^d}\int_{\mathbb{T}}\int_{\mathbb{R}}\int_{\mathbb{R}^d}J(z,r)\mu(\xi,q)\mu(\xi-z,q-r)\nonumber\\
&&\cdot\Big|\nabla u(x,t)\Big(-z+\chi_1(\xi-z,q-r)-\chi_1(\xi,q)\Big)\Big|^{p-2}\nonumber\\
&&\cdot\Big[-z\otimes\chi_1(\xi-z,q-r)+{\frac{1}{2}}z\otimes z\Big]\nabla\nabla u(x,t)dzdrd\xi dq,
\end{eqnarray}
\begin{eqnarray}\label{Equ2.24}
\lim_{\varepsilon\rightarrow 0}I_2^{\varepsilon}&=&-(p-1)\frac{\partial u(x,t)}{\partial t}\int_{\mathbb{T}^d}\int_{\mathbb{T}}\int_{\mathbb{R}}\int_{\mathbb{R}^d}J(z,r)\mu(\xi,q)\mu(\xi-z,q-r)\nonumber\\
&&\cdot\Big|\nabla u(x,t)(-z+\chi_1(\xi-z,q-r)-\chi_1(\xi,q))\Big|^{p-2}rdzdrd\xi dq.
\end{eqnarray}
For \cite[Theorem 24,Lemma 23]{Rossi2017} there exists functions $u_{\varepsilon_{k}}$ solutions of the problem $ \mathbb{S}(J_{\varepsilon_k},\nu,w^{\varepsilon},f)$ that converges uniformly to a function $u$ when $k$ goes to infinity.

Since  $u$ is a viscosity solution of the local problem. Let $\varphi(x, t)\in C^{2}\left(\mathbb{R}^{d+1}\right)$, and assume that $\bar{u}-\varphi$ has an strict maximum at the point $(x, t) \in \mathbb{R}^{d} \times(0, \infty)$ and $\nabla \varphi(x, t) \neq 0 .$ Then, as $u_{\varepsilon_{k}} \rightarrow u$ uniformly in a neighborhood of $(x, t)$, there exist $\left(x_{\varepsilon_{k}}, t_{\varepsilon_{k}}\right) \in \mathbb{R}^{d} \times(0, \infty)$
such that maximum point (x,t) $\to  \overline{u_{\varepsilon_{k}}}-(\varphi+\varepsilon\chi_1(\frac{x}{\varepsilon},\frac{t}{\varepsilon^2})\nabla \varphi(x,t)
+\varepsilon^2\chi_2(\frac{x}{\varepsilon},\frac{t}{\varepsilon^2})\nabla \nabla \varphi(x,t))$ in the set  $ B^\varepsilon=\left (B_{R_{{\varepsilon_{k}}}}\left(x_{\varepsilon_{k}}\right)  , ( t_{\varepsilon_{k}}-R_{{\varepsilon_{k}}}, t_{\varepsilon_{k}}+R_{{\varepsilon_{k}}}) \cap(0,\infty)\right)$, $ R_{{\varepsilon_{k}}}\to \infty$  with $\left(x_{\varepsilon_{k}}, t_{\varepsilon_{k}}\right) \rightarrow(x, t)$ as $\varepsilon_{k} \to 0$, where $ \chi_1,\chi_2 \in C_b(\mathbb{T}^d\times \mathbb{T})$ are introduced in lemmata 2.2 and  lemmata 2.3. Then, we denote test function $\varphi^\varepsilon= \varphi+\varepsilon\chi_1(\frac{x}{\varepsilon},\frac{t}{\varepsilon^2})\nabla \varphi(x,t)
+\varepsilon^2\chi_2(\frac{x}{\varepsilon},\frac{t}{\varepsilon^2})\nabla \nabla \varphi(x,t)$ and have
$$
\overline{u_{\varepsilon_{k}}}(y, s)-u_{\varepsilon_{k}}\left(x_{\varepsilon_{k}}, t_{\varepsilon_{k}}\right) \leq \varphi^\varepsilon(y, s)-\varphi^\varepsilon\left(x_{\varepsilon_{k}}, t_{\varepsilon_{k}}\right) \quad as~ \varepsilon <<1.
$$
Hence,
$$
F\left(\overline{u_{\varepsilon_{k}}}(y, s)-u_{\varepsilon_{k}}\left(x_{\varepsilon_{k}}, t_{\varepsilon_{k}}\right)\right) \leq F\left(\varphi^\varepsilon(y, s)-\varphi^\varepsilon\left(x_{\varepsilon_{k}}, t_{\varepsilon_{k}}\right)\right)
$$
This implies that

\begin{eqnarray*}
	0 &=&\frac{1}{\varepsilon_{k}^{p}} \int_{\mathbb{R}^{d} \times \mathbb{R}}  J_{\varepsilon_{k}}\left(\left(x_{\varepsilon_{k}}-y\right),\left(t_{\varepsilon_{k}}-s\right)\right)\nu_{\varepsilon_{k}}(x_{\varepsilon_{k}},y,t_{\varepsilon_{k}},s) F\left(\overline{u_{\varepsilon_{k}}}(y, s)-u_{\varepsilon_{k}}\left(x_{\varepsilon_{k}}, t_{\varepsilon_{k}}\right)\right) d y d s \nonumber\\
	& \leq& \frac{1}{\varepsilon_{k}^{p}}\int_{\mathbb{R}^{d} \times \mathbb{R}} J(z, r) \nu_{\varepsilon_{k}}(x_{\varepsilon_{k}},x_{\varepsilon_{k}}-\varepsilon_{k},t_{\varepsilon_{k}}, t_{\varepsilon_{k}}-\varepsilon_{k}^{2} r) F\left(\varphi^{\varepsilon_k}\left(x_{\varepsilon_{k}}-\varepsilon_{k} z, t_{\varepsilon_{k}}-\varepsilon_{k}^{2} r\right)-\varphi^{\varepsilon_k}\left(x_{\varepsilon_{k}}, t_{\varepsilon_{k}}\right)\right) d z d r,
\end{eqnarray*}
we hope the right term of the last inequality will converge to a  local equation which acting over $\varphi^{\varepsilon_k}$ when $\varepsilon_{k}$ tends to 0.
Therefore from (2.23) and (2.24),  we directly get 

$$
0\leq-\mathfrak{N}(\nabla \varphi)\frac{\partial \varphi}{\partial t}+\mathfrak{P}(\nabla \varphi,\nabla \nabla \varphi) .
$$

On the other hand if $\psi(x, t)$ is a $C^{2}\left(\mathbb{R}^{d+1}\right)$ function such that $\bar{u}-\psi$ has a strict minimum at $(x, t) \in$ $\mathbb{R}^{d} \times(0, \infty)$ and $\nabla \psi(x, t) \neq 0$, arguing as before we obtain
$$
0 \geq \frac{1}{\varepsilon_{k}^{p}} \int_{\mathbb{R}^{d} \times \mathbb{R}} J(z, r)  \nu(\frac{x_{\varepsilon_{k}}}{\varepsilon_k},\frac{x_{\varepsilon_{k}}}{\varepsilon_{k}}-z,\frac{t_{\varepsilon_{k}}}{\varepsilon^2_{k}},\frac{t_{\varepsilon_{k}}}{\varepsilon^2_{k}}-r)F\left(\psi^{\varepsilon_k}\left(x_{\varepsilon_{k}}-\varepsilon_{k} z, t_{\varepsilon_{k}}-\varepsilon_{k}^{2} r\right)-\psi^{\varepsilon_k}(x_{\varepsilon_k}, t_{\varepsilon_k})\right) d z d r.
$$

Using (2.23) and (2.24) again, we get
\begin{eqnarray}\label{Equ2.25}
	0\geq-\mathfrak{N}(\nabla \psi)\frac{\partial \psi}{\partial t}+\mathfrak{P}(\nabla \psi,\nabla \nabla \psi).
\end{eqnarray}
It's easy to check that $w^{\varepsilon} \to u$ in the sense of viscosity solution  as $\varepsilon\rightarrow 0^{+},$ and $\lim\limits_{\varepsilon \rightarrow 0} w^{\varepsilon}(x,0)=f(x)$.  From above we conclude that $u$ is a viscosity solution to the effective problem
\begin{eqnarray}\label{Equ2.26}
	\left\{\begin{array}{ll}
		\mathfrak{N}(\nabla u)\frac{\partial u}{\partial t}=\mathfrak{P}(\nabla u,\nabla\nabla u),\; &(x,t)\in\mathbb{R}^d\times (0,T), \\
		u(x,0)=f(x),&x\in \mathbb{R}^d. \end{array}\right.
\end{eqnarray}
This ends the proof of Theorem \ref{The2.1}. $\Box$
\begin{remark}\label{rem2.1} In one dimensional case $d=1$, (\ref{Equ2.26}) becomes
\begin{eqnarray}\label{Equ2.27}
\left\{\begin{array}{ll}
\mathfrak{N}|\frac{\partial u}{\partial x}|^{p-2}\frac{\partial u}{\partial t}=\mathfrak{P}\frac{\partial}{\partial x}(|\frac{\partial u}{\partial x}|^{p-2}\frac{\partial u}{\partial x}),\;&(x,t)\in\mathbb{R}\times[0,\infty), \\
u(x,0)=f(x), &x\in\mathbb{R}, \end{array}\right.
\end{eqnarray}
where $\mathfrak{N}$ and $\mathfrak{P}$ are constants depending on $J$ that are given by
\begin{eqnarray*}
\mathfrak{N}&=&(p-1)\int_{\mathbb{T}}\int_{\mathbb{T}}\int_{\mathbb{R}}\int_{\mathbb{R}}J(z,r)\mu(\xi,q)\mu(\xi-z,q-r)\nonumber\\
&&\cdot|\chi_1(\xi-z,q-r)-\chi_1(\xi,q)-z|^{p-2}rdzdrd\xi dq, \nonumber\\
\mathfrak{P}&=&\int_{\mathbb{T}}\int_{\mathbb{T}}\int_{\mathbb{R}}\int_{\mathbb{R}}J(z,r)\mu(\xi,q)\mu(\xi-z,q-r)|\chi_1(\xi-z,q-r)-\chi_1(\xi,q)-z|^{p-2}\nonumber\\
&&\cdot[-z\otimes\chi_1(\xi-z,q-r)+{\frac{1}{2}}z\otimes z]dzdrd\xi dq.
\end{eqnarray*}
\end{remark}
\begin{remark}\label{rem2.2} When $\mu\equiv1$, it's easy to check that $\chi_1=\chi_2=0.$ It degenerates to the initial value problem given in \cite[Theorem 3]{Rossi2017}
\begin{eqnarray}\label{Equ2.28}
\left\{\begin{array}{l}
\mathfrak{N}|\nabla u|^{p-2}\frac{\partial u}{\partial t}=\mathfrak{P}(|\nabla u|^{p-2}\Delta u+(p-2)|\nabla u|^{p-4}\nabla u D^2u\nabla u), \\
u(x,0)=f(x), \end{array}\right.
\end{eqnarray}
for $(x,t)\in \mathbb{R}^d\times(0, +\infty)$, where $\mathfrak{N}$ and $\mathfrak{P}$ depending on $J$ that are given by
\begin{eqnarray*}
\mathfrak{N}=\int_{\mathbb{R}}\int_{\mathbb{R}^d}J(z,r)|z_1|^{p-2}rd zdr, \; \mathfrak{P}=\frac{1}{2}\int_{\mathbb{R}}\int_{\mathbb{R}^d}J(z,r)|z_1|^{p}dzdr,
\end{eqnarray*}
and $z_1$ is the first component of $z=(z_1,z_2,\cdots,z_d)$.
\end{remark}

\section{$\mu(x,t)$ is periodic in $x$ and stationary in $t$}\label{sec3}
\setcounter{equation}{0}

Assume that $\nu(\frac{x}{\varepsilon},\frac{t}{\varepsilon^2};\frac{y}{\varepsilon},\frac{s}{\varepsilon^2})=\mu(\frac{x}{\varepsilon},\frac{t}{\varepsilon^2})\mu(\frac{y}{\varepsilon},\frac{s}{\varepsilon^2})$ is periodic in $x,y$ and stationary in $t,s$. We consider the nonlocal linear problem (\ref{Equ2.1}): $\mathbb{S}(J_{\varepsilon},\nu,u^{\varepsilon},f)=0$ for $p=2$, that is,
\begin{eqnarray}\label{Equ3.1}
\frac{1}{\varepsilon^{d+4}}\int\int_{\mathbb{R}^{d}\times\mathbb{R}}J(\frac{x-y}{\varepsilon},\frac{t-s}{\varepsilon^2})\mu(\frac{x}{\varepsilon},\frac{t}{\varepsilon^2}) \mu(\frac{y}{\varepsilon},\frac{s}{\varepsilon^2})(\bar{u}^{\varepsilon}(y,s)-u^{\varepsilon}(x,t))dyds=0,
\end{eqnarray}
where $\bar{u}^{\varepsilon}$ is defined in (\ref{Equ1.2}). For the classical initial value problem, we have $\bar{u}(y,s)=u(y,s)$. Assume that the initial function $f$ is at least four times continuously differentiable and $u_0(x,t)$ has enough regularity in $x$, that is, for any $l>0$ there is a constant $C_{l}(T)>0$ such that a solution of problem (\ref{Equ3.1}) satisfies the estimate
\begin{eqnarray}\label{Equ3.2}
\sum_{|\mathbf{k}|\leq 4}|\partial_{x}^{\mathbf{k}}u_0(x,t)|\leq C_{l}(T)(1+|x|)^{-l}, \forall\; (x,t)\in\mathbb{R}^{d}\times[0, T].
\end{eqnarray}

For the sake of employing invariant principle to some stochastic process, we introduce the so-called maximum correlation coefficient.  Let $(\Omega, \mathcal{F}, \mathbf{P})$ be a standard probability space equipped with a measure preserving ergodic
dynamical system  $T_t, t\in \mathbb{R}$. Setting $\mathcal{F}_{\leq s}=\sigma\{\mu(x,t): t\leq s\}$ and $\mathcal{F}_{\geq s}=\sigma\{\mu(x,t): t\geq s\}$. We define the maximum correlation coefficient
\begin{eqnarray}\label{Equ3.3}
\rho(s)=\sup_{\varsigma_{1},\varsigma_{2}}\mathbf{E}(\varsigma_{1}\varsigma_{2}),
\end{eqnarray}
where the supremum is taken over all $\mathcal{F}_{\leq s}$-measurable $\varsigma_{1}$ and $\mathcal{F}_{\geq s}$-measurable $\varsigma_{2}$ such that $\mathbf{E}\varsigma_{1}=\mathbf{E}\varsigma_{2}=0$, and $\mathbf{E}\varsigma_{1}^{2}=\mathbf{E}\varsigma_{2}^{2}=1$.

Throughout this paper we assume that the function $\rho$ satisfies $\int_{0}^{\infty}\rho(s)ds<+\infty$. Denote by $\mathcal{F}_{\leq T}^{\mu,h}$ the $\sigma-$algebra $\sigma\{\mu(z,s),h(z,s):s\leq T\}$. The $\sigma-$algebra $\mathcal{F}_{\geq T}^{\mu,h}$ is defined accordingly. Let $\rho_{\mu,h}(r)$ be maximum correlation coefficient of $(\mu,h)$. Denote
\begin{eqnarray}\label{Equ3.4}
\digamma(q)&=&\int_{\mathbb{T}^d}\Big[\int_{\mathbb{R}}\mu(\xi,q)\int_{\mathbb{R}^d}J(z,r)\mu(\xi-z,q-r)z\otimes \chi(\xi-z,q-r)dzdr\nonumber\\
&&-\mathbf{E}\int_{\mathbb{R}}\mu(\xi,q)\int_{\mathbb{R}^d}J(z,r)\mu(\xi-z,q-r)z\otimes \chi(\xi-z,q-r)dzdr\Big]d\xi.
\end{eqnarray}
\begin{remark}\label{rem3.1} The condition $\int_{0}^{\infty}\rho(s)ds<+\infty$ on maximum correlation coefficient (\ref{Equ3.3}) may not be very intuitive. In fact, we consider a diffusion process $\eta_{s}, s\in(-\infty,+\infty)$, with values in $\mathbb{R}^{d}$ or on a compact manifold. This process is defined on a probability space $(\Omega, \mathcal{F}, \mathbb{P})$. The corresponding Ito's equation is written as
\begin{eqnarray}\label{Equ3.5}
d\eta_{t}=b(\eta_{t})dt+\sigma(\eta_{t})dB_{t},
\end{eqnarray}
where $B_{t}$ stands for a standard $d-$dimensional Wiener process, and the process $\eta$ has a unique invariant probability measure when the infinitesimal generator of $\eta$ satisfies some conditions. More details can be seen in \cite{Pardoux}.
\end{remark}
Denote
\begin{eqnarray}
&&\mathfrak{L}\chi(\xi,q)\equiv\int_{\mathbb{R}}\int_{\mathbb{R}^d}J(z,r)\mu(x-z,t-r)\Big(\chi(\xi-z,q-r)-\chi(\xi,q)\Big)dzdr, \label{Equ3.6}\\
&&h(\xi,q)=\int_{\mathbb{R}}\int_{\mathbb{R}^d}J(z,r)z\mu(\xi-z,q-r)dzdr.\label{Equ3.7}
\end{eqnarray}
\begin{lemma}\label{lem3.1} For some $\varrho,C>0$ and the vector-valued function $\digamma(\cdot)$ we have
\begin{eqnarray}\label{Equ3.8}
\left\|\mathbb{E}\left\{\digamma(s)\mid\mathcal{F}_{\leq 0}^{\mu,h}\right\}\right\|_{L^{2}(\Omega)}\leq C(\rho_{\mu,h}(\frac{s}{2})+e^{-\frac{\varrho s}{2}}).
\end{eqnarray}
\end{lemma}
\noindent{\bf Proof.} Without loss of generality, we give a specific kernel function $J$. For example, the Weierstrass kernel $\mathcal{W}(x,t)=(4\pi t)^{-d/2}e^{-\frac{|x|^{2}}{4t}}$ for $x\in\mathbb{R}^{d}$ and $t>0$, set
\begin{eqnarray}\label{Equ3.9}
\mathbf{I}=\left\{(x,t)\in \mathbb{R}_{+}^{d+1}: \mathcal{W}(x,t)\geq 1\right\}\;{and}\; J(x,t)=\frac{1}{4}\mathbb{I}_{\mathbf{I}}(x,t)\frac{|x|^{2}}{t^{2}},
\end{eqnarray}
where $\mathbb{I}_{\mathbf{I}}(x,t)$ is an indicator function on the set $\mathbf{I}$. We represent the Green function $\iota(z,s)$ of $\mathfrak{L}$ defined in (\ref{Equ3.6}) on the interval $0\leqslant s\leqslant\hat{s}$ as a sum $\iota(z,s)=\iota^{1}(z,s)+\iota^{2}(z,s)$, where $\iota^{1}$ and $\iota^{2}$ satisfy that
$$
\mathfrak{L}\iota^{i}(z, s)=0,\; s<\hat{s},\left.\iota^{1}\right|_{s=\hat{s}}=1,\left.\iota^{2}\right|_{s=\hat{s}}=\iota(z,\hat{s})-1, i=1,2.
$$
From the definition of $J$ in (\ref{Equ3.9}) and \cite{Rossi2015}, it follows that for all $\hat{s}\geq t_0+1$,
$$
 (\mathfrak{P}) \quad\quad\|\iota(\cdot,x_0,\cdot,t_0)-1\|_{L^2\left((\hat{s},\hat{s}+1); L^2\left(\mathbb{T}^d\right)\right)}\leq Ce^{-\varrho(\hat{s}-t_0)},
$$
we proved the inequality in the appendix. Choosing $\chi(z,a)=\chi^{1}(z,a)+\chi^{2}(z,a)$ such that
\begin{eqnarray}\label{Equ3.10}
\left\{\begin{array}{l}
\chi^{(1)}(z,a)=\displaystyle{\int_{-\infty}^{a-1}\int_{\mathbb{T}^{d}}(\iota(z,z_1,s,s_1)-1)h(z_1,s_1)dz_1ds_1,}\\
\chi^{(2)}(z,a)=\displaystyle{\int_{a-1}^{a}\int_{\mathbb{T}^{d}}(\iota(z,z_1,s,s_1))h(z_1,s_1)dz_1ds_1,}  \end{array}\right.
\end{eqnarray}
and putting (\ref{Equ3.10}) into (\ref{Equ3.4}), for $i=1,2$, we have
\begin{eqnarray}\label{Equ3.11}
\digamma^{(i)}(q)&=&\int_{\mathbb{T}^d}\Big[\int_{\mathbb{R}}\mu(\xi,q)\int_{\mathbb{R}^d}J(z,r)\mu(\xi-z,q-r)z\otimes \chi^{(i)}(\xi-z,q-r)dz dr\nonumber\\
&-&\mathbf{E}\int_{\mathbb{R}}\mu(\xi,q)\int_{\mathbb{R}^d}J(z,r)\mu(\xi-z,q-r)z\otimes\chi^{(i)}(\xi-z,q-r)dzdr\Big]d\xi.
\end{eqnarray}
Using the maximum principle in \cite[Theorem 3]{Aimar2018} and similar to the proof of \cite[Lemma 3]{Pyatnitskii}, we can obtain the desired inequality (\ref{Equ3.8}). $\Box$
\begin{lemma}\label{lem3.2} If $h\in L^2_{loc}(\mathbb{R}, L^2(\mathbb{T}^d))$ satisfies $\displaystyle{\int_{\mathbb{T}^{d}}h(\xi,q)d\xi=0~a.s.}$, then the equation
\begin{eqnarray}\label{Equ3.12}
\mathfrak{L}\chi(\xi,q)=h,\;\xi\in\mathbb{T}^d, q\in \mathbb{R},
\end{eqnarray}
has a stationary solution $\chi\in L_{\mathrm{loc}}^{2}(\mathbb{R}, L^{2}(\mathbb{T}^{d}))^d$.
Furthermore, if $h\in L^\infty(\mathbb{R}\times \mathbb{T}^d)$, then we have $\|\chi\|_{L^{\infty}(\mathbb{R}\times\mathbb{T}^{d})}\leq C,$ where $C$ is a non-random positive constant.
\end{lemma}
{\bf Proof.} First we consider the following problem:
\begin{eqnarray*}
\left\{\begin{array}{l}
\mathfrak{L}\chi^T(\xi,q)=h,\;(\xi,q)\in \mathbb{ T}^d\times(-\infty,T],\\ \chi^T(\xi,q)|_{q=T}=\chi_0(\xi),  \end{array}\right.
\end{eqnarray*}
where $T \in \mathbb{R}$ and $\chi_0 \in L^2(\mathbb{T}^d)$ satisfies the equation $\displaystyle{\int_{\mathbb{T}^d}\chi_0(\xi)d\xi=0.}$

The existence of solution $\chi^T-\chi^{T-1}$ can be obtained from \cite{Piatnitski2017}. The key idea of proof is to decompose the operator $\mathfrak{L}$ in (\ref{Equ3.6}) into the sum of a positive invertible operator and a compact operator, and then one can apply the Fredholm theorem to get solution $\chi(\xi,q)$, it is periodic in $\xi$ and stationary in $q$. The uniqueness of solution $\chi(\xi,q)$ for local parabolic equation can be obtained following the proof of \cite[Lemma 4]{Pyatnitskii}, since $\chi^T$ converges exponentially to $\chi$ in (\ref{Equ3.12}) (can be seen in Section 6) as $T\to +\infty$, $\chi(\xi,q)$ is also unique up to an additive (random) constant. $\Box$

As complexity of calculations we just consider the linear equation (\ref{Equ3.1}) with $p=2$. For $p\neq2$, we will find a new method and consider in the forth coming paper.

Similar to Theorem \ref{The2.1}, if $J$ satisfies the condition (\ref{Equ1.3}), is compactly supported in $\mathbb{R}^{d+1}$ and it is radial in $x\in \mathbb{R}^d$ for each $t$, we easily get the following qualitative result.
\begin{theorem}\label{The3.1} If a bounded function $f\in C^{2}(\mathbb{R}^{d},\mathbb{R})$ has bounded derivatives, the solution of (\ref{Equ3.1}) converges in $L_{\omega}^{2}((0,T)\times\mathbb{R}^d)$ as $\varepsilon\rightarrow 0^{+}$ to a solution of the initial value problem
\begin{eqnarray}\label{Equ3.13}
\left\{\begin{array}{l}
\hat{\alpha}\frac{\partial u_0}{\partial t}+\Theta\cdot\nabla\nabla u_0(x,t)=0, \\ u_0(x,0)=f(x),  \end{array}\right.
\end{eqnarray}
where $\hat{\alpha}$ and $\Theta$ depending on $J$ that are given by
\begin{eqnarray}
&&\hat{\alpha}=-\mathbb{E}\int_{\mathbb{T}^{d}}\int_{\mathbb{R}}\mu(\xi,q)\int_{\mathbb{R}^{d}}J(z,r)\mu(\xi-z,q-r)rdzdrd\xi,\label{Equ3.14}\\
&&\Theta=\mathbb{E}\int_{\mathbb{T}^{d}}\int_{\mathbb{R}}\mu(\xi,q)\int_{\mathbb{R}^{d}}J(z,r)\mu(\xi-z,q-r)\nonumber\\
&&\qquad\cdot\Big(-z\otimes\chi(\xi-z,q-r)+{\frac{1}{2}z\otimes z}\Big)dzdrd\xi, \label{Equ3.15}
\end{eqnarray}
and $\Theta\cdot\nabla\nabla u_0=\sum\limits_{i,j=1}^{d}\Theta_{ij}\partial_{x_{i}x_{j}}u_0(x,t)$, $\chi$ will be introduced later.
\end{theorem}

If we simplify $\mu$, we can calculate the determined coefficients $\hat{\alpha},\Theta$. Assume that random statistically homogeneous ergodic media $\mu_\varepsilon=\mu(\frac{t}{\varepsilon^2},\omega), t \in [0,1]$ and periodization in 1-dimensional media  which is
depicted in Figure 1, where $\eta_1=\frac{1}{2}$ and $\eta_2=2$ and with probability $p=\frac{1}{2}.$ More examples can be seen in \cite{periodization}.

Next we can directly compute
\begin{eqnarray*}
\Theta&=&\mathbb{E}\int_{\mathbb{T}}\int_{\mathbb{R}}\mu(q)\int_{\mathbb{R}^{d}}J(z,r)\mu(q-r)\Big(-z\chi(q-r)+{\frac{1}{2}z^2}\Big)dzdrd\xi\nonumber\\
&=& \int_{\mathbb{R}}\mathbb{E}\mu(0,\omega)\frac{1}{2}\int_{\mathbb{R}}\mu(-r,\omega)J(z,r){z^2}dzdr.
\end{eqnarray*}
If we choose kernel function $J$ in (\ref{Equ3.9}),  when $d$=1 and $\mu=1$ after simple calculations get
\begin{eqnarray}\label{Equ3.16}
\mathbf{k}_0=\int\int J(z,r)rdzdr=\frac{1}{2}\int\int J(z,r){z^2}dzdr=\frac{1}{36 \sqrt{3}\pi},
\end{eqnarray}
hence we have $\displaystyle{\hat{\alpha}=-\int_{\mathbb{R}}\int_{\mathbb{R}}J(z,r)rdzdr=-\mathbf{k}_0=-\Theta.}$
\begin{figure}[htb]
  \centering
  \includegraphics[width=0.7\textwidth]{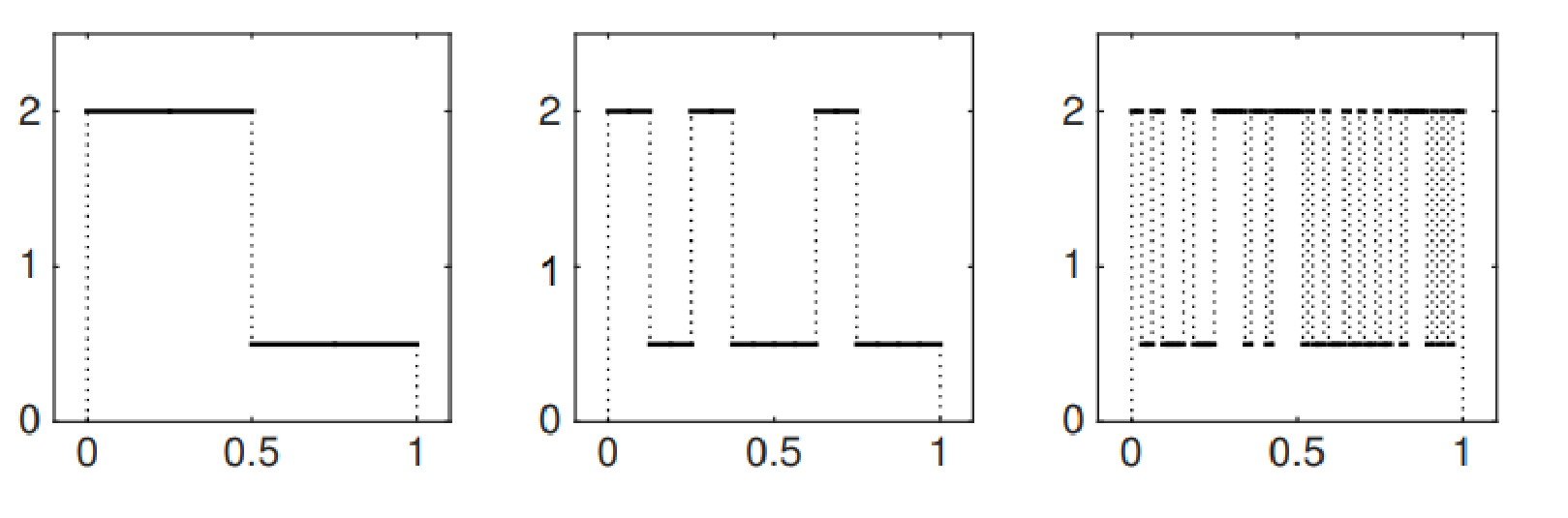}
  \caption{{\small One realization of the one-dimensional tile-based random media w.r.t. tile size $\varepsilon =2^{-1},2^{-2},2^{-3}$.}}
  \label{Fig.1}
\end{figure}

For the mesoscopic form from non-locality to locality, we now consider the following case
\begin{eqnarray*}
\frac{1}{\varepsilon^{2}}\int_{\mathbb{R}}\int_{\mathbb{R}}\int_{\mathbb{R}^{d}}\int_{\mathbb{R}^{d}}J(z,r)\Big(v^{\varepsilon}(x-\varepsilon z,t-\varepsilon^2 r)-v^{\varepsilon}(x,t)\Big)^{2}dzdxdtdr\leq C.
\end{eqnarray*}
For convenience we assume that $Q_T=[-1,1]^{d+1}$ (sphere and cube are equivalent), $\phi(s)\in C_{0}^{\infty}(\mathbb{R}^d)$ function such that $0\leq\phi\leq 1$, $\phi(s)=1$ for $s\in B(0,1)$, $\phi(s)=0$ for $s \in B^{c}(0,2)$, and $|\phi^{\prime}(s)|\leq 2$, $\eta (t) \in C_{0}^{\infty}(0,T)$. Denote $\bar{v}^{\varepsilon}(x,t)=\phi(x)\eta (t)v^{\varepsilon}(x,t)$. It is easy to check that
\begin{eqnarray*}
\frac{1}{\varepsilon^{2}}\int_{\mathbb{R}}\int_{\mathbb{R}}\int_{\mathbb{R}^{d}}\int_{\mathbb{R}^{d}} J(z,r)\Big(\bar{v}^{\varepsilon}(x-\varepsilon z,t-\varepsilon^{2}r)-\bar{v}^{\varepsilon}(x,t)\Big)^{2}dzdxdtdr\leq C.
\end{eqnarray*}
Denote $\tilde{Q}=[-\frac{1}{2},\frac{1}{2}]^{d+1}$, $\hat{v}^{\varepsilon}(x,t)$ is the $\tilde{Q}-$periodic extension of $\bar{v}^{\varepsilon}(x,t)$, and
\begin{eqnarray*}
\mathcal{P}=\{(x,t)\in \mathbb{R}_{+}^{d+1}:J(x,t)\geq \frac{1}{2}J_{max}\wedge 1\},
\end{eqnarray*}
where $J_{max}=\sup_{(x,t)\in\mathbb{R}^d\times(0,T)}J(x,t)$ and $\tilde{J}(z,r)=1_{\{(z,r)\in \mathcal{P}\}}J(z,r)(\int_{\mathcal{P}}J(z,r)dzdr)^{-1}.$ For the extended function we have
\begin{eqnarray*}
\frac{1}{\varepsilon^{2}}\int_{\tilde{Q}}\int_{\mathbb{R}^{d}} \tilde{J}(z,r)\Big(\hat{v}^{\varepsilon}(x-\varepsilon z,s-\varepsilon^2 r)-\hat{v}^{\varepsilon}(x,s)\Big)^{2}dzdrdxds\leq C.
\end{eqnarray*}
The functions $e_{k,l}(x,t)=e^{i2\pi(kx+lt)}(k\in Z^{d},l\in Z)$ form an orthonormal basis of $L^{2}(Q_T)$, and
\begin{eqnarray*}
&&\hat{v}^{\varepsilon}(x,t)=\sum_{k,l}\alpha_{k,l}^{\varepsilon}e_{k,l}(x,t),\;\hat{v}^{\varepsilon}(x-\varepsilon z,t-\varepsilon^2 r)=\sum_{k,l}\alpha_{k,l}^{\varepsilon}e^{-i2\pi(\varepsilon kz+\varepsilon^2 lr)}e_{k,l}(x,t), \nonumber\\
&&\|\hat{v}^{\varepsilon}(x,t)\|^{2}=\sum_{k,l}(\alpha_{k,l}^{\varepsilon})^{2}, \nonumber\\
&&\|\hat{v}^{\varepsilon}(x-\varepsilon z,t-\varepsilon^2 r)-\hat{v}^{\varepsilon}(x,t)\|^{2}=\sum_{k,l}(\alpha_{k,l}^{\varepsilon})^{2}|e^{-i2\pi(\varepsilon kz+\varepsilon^2 lr)}-1|^{2}.
\end{eqnarray*}
Then we obtain
\begin{eqnarray*}
\frac{1}{\varepsilon^{2}}\sum_{k,l}(\alpha_{k,l}^{\varepsilon})^{2}\int_{\mathbb{R}}\int_{\mathbb{R}^{d}}\tilde{J}(z,r)|e^{-i2\pi(\varepsilon kz+\varepsilon^2 lr)}-1|^{2}dzdr\leq C.
\end{eqnarray*}
\begin{lemma}\label{lem3.3} For any $k\in Z^{d},l\in Z $ and any $0<\varepsilon<1$ there exist constants $C_{1}$, $C_{2}$ (depending on $d$) such that
\begin{eqnarray*}
W=\Big\{\hat{v}^{\varepsilon}(x,t)\Big|\frac{1}{\varepsilon^{2}}\sum_{k,l}\Big(\alpha_{k,l}^{\varepsilon}\Big)^{2}\int_{\mathbb{R}}\int_{\mathbb{R}^{d}}\tilde{J}(z,r)\Big|e^{-i 2\pi(\varepsilon kz+\varepsilon^2 lr)}-1\Big|^{2}dzdr\leq C\Big\}
\end{eqnarray*}
is compact in $L^2(Q_T)$.
\end{lemma}
\begin{remark}\label{rem3.2} The proof is similar in \cite{Piatnitski2022} and we don't repeat it. This lemma implies an important idea, that is, the low-frequency part of the Fourier form of the solution of the non-local operator equation can control the high-frequency domain. We can construct a finite-dimensional sub-space which has finite $\delta-$set cover high-frequency part, so the low-frequency part has the differentiable regularity that the local operator has , and this is precisely the critical connection between the non-local operator and the local operator. A. Piatnitski et al \cite{Piatnitski2022} studied the convergence rate of nonlocal elliptic operators in by using the spectral methods, the proof of convergence in Lemma \ref{lem3.5} will be used it later.
\end{remark}

We next deal with the formal asymptotic expansion of solution $u^{\varepsilon}(x,t)$ to the problem (\ref{Equ3.1}) with $p=2$, inspired by the paper \cite{Piatnitski2015} our goal is to obtain the structure of the leading terms of the difference $u^{\varepsilon}(x,t)-u_0(x,t)$, which is the main result in this section. In order to use the multi-scale asymptotic expansion method we consider $\xi=x/\varepsilon$ and $q=t/\varepsilon^{2}$ as independent variables and have formulas
\begin{eqnarray*}
&&\partial_x f(x,\frac{x}{\varepsilon})=\partial_x f(x,\xi)+\frac{1}{\varepsilon}\partial_{\xi} f(x,\xi),\;\xi=\frac{x}{\varepsilon},\\
&&\partial_t g(t,\frac{t}{\varepsilon^2})=\partial_q g(t,q)+\frac{1}{\varepsilon^2}\partial_q g(t,q),\;q=\frac{t}{\varepsilon^2}.
\end{eqnarray*}

We represent a solution $u^{\varepsilon}(x,t)$ as the following asymptotic series in integer powers of $\varepsilon$ by using the classical multi-scale method
\begin{eqnarray}\label{Equ3.17}
u^{\varepsilon}(x,t)=u_{0}(x,t)+\varepsilon u_{1}(x,t,\frac{x}{\varepsilon},\frac{t}{\varepsilon^{2}})+\varepsilon^{2}u_{2}(x,t,\frac{x}{\varepsilon},\frac{t}{\varepsilon^{2}})+\cdots,
\end{eqnarray}
where all the functions $u_{j}(x,t,\xi,q)$ are periodic in $\xi$ but not always stationary in $q$. For brevity of notation $u_0$ sometimes be written as $u$. Substituting the expression on the left-hand side of (\ref{Equ3.1}) for $u^{\varepsilon}$ in (\ref{Equ3.17}) and collecting power-like terms in $\varepsilon^{-1}, \varepsilon^{0}, \cdots$ yield
\begin{eqnarray}\label{Equ3.18}
(\varepsilon^{-1}):\;T_{-1}&\equiv&\int_{\mathbb{R}}\mu(\frac{x}{\varepsilon},\frac{t}{\varepsilon^2})\int_{\mathbb{R}^d}J(z,r)\mu(\frac{x}{\varepsilon}-z,\frac{t}{\varepsilon^2}-r)[-z\nabla u(x,t)\nonumber\\
&+&u_1(x,t,\frac{x}{\varepsilon}-z,\frac{t}{\varepsilon^2}-r)-u_1(x,t,\frac{x}{\varepsilon},\frac{t}{\varepsilon^2})]dzdr,
\end{eqnarray}
\begin{eqnarray}\label{Equ3.19}
(\varepsilon^{0}):\; T_{0}&\equiv&\int_{\mathbb{R}}\mu(\frac{x}{\varepsilon},\frac{t}{\varepsilon^2})\int_{\mathbb{R}^d}J(z,r)\mu(\frac{x}{\varepsilon}-z,\frac{t}{\varepsilon^2}-r)\nonumber\\
&&\cdot[-z\otimes \partial_x u_1(x,t,\frac{x}{\varepsilon}-z,\frac{t}{\varepsilon^2}-r)-u_2(x,t,\frac{x}{\varepsilon},\frac{t}{\varepsilon^2})\nonumber\\
&&+u_2(x,t,\frac{x}{\varepsilon}-z,\frac{t}{\varepsilon^2}-r)+{\frac{1}{2}}z^{2}\cdot\nabla\nabla u-r\frac{\partial u(x,t)}{\partial t}]dzdr.
\end{eqnarray}

Our aim is to prove that $T_{-1}=0$ and $\lim\limits_{\varepsilon\rightarrow 0}T_{0}$ exists and satisfies a limit equation. By Lemma \ref{lem3.2} we know that the equation $T_{-1}=0$ has a unique stationary solution $u_1$ such that $u_{1}(x,t,\xi,q)=\chi(\xi,q)\nabla u(x,t),$ where $\chi=\{\chi^{j}(\xi,q)\}_{j=1}^{d}$ is a vector-valued function and satisfies that
\begin{eqnarray}\label{Equ3.20}
\int_{\mathbb{R}}\int_{\mathbb{R}^d}J(z,r)\mu(\xi-z,q-r)(-z+\chi(\xi-z,q-r)-\chi(\xi,q))dzdr=0.
\end{eqnarray}
By Lemma \ref{lem3.2} we have $\chi\in L_{\mathrm{loc}}^{2}(\mathbb{R},L^{2}(\mathbb{T}^{d}))^d$ and $\int_{\mathbb{T}^d}\chi(z,s)dz=0$. Denote
\begin{eqnarray}\label{Equ3.21}
\Pi_1(q)&=&\int_{\mathbb{T}^d}[\int_{\mathbb{R}}\mu(\xi,q)\int_{\mathbb{R}^d}J(z,r)\mu(\xi-z,q-r)\nonumber\\
&&\cdot(-z\otimes\chi(\xi-z,q-r)+{\frac{1}{2}}z^{2})dzdr-\Theta]d\xi,
\end{eqnarray}
\begin{eqnarray}\label{Equ3.22}
\Pi_2(\xi,q)&=&[\int_{\mathbb{R}}\mu(\xi,q)\int_{\mathbb{R}^d}J(z,r)\mu(\xi-z,q-r)\nonumber\\
&&\cdot(-z\otimes\chi(\xi-z,q-r)+{\frac{1}{2}}z^{2})dzdr-\Theta]-\Pi_1,
\end{eqnarray}
\begin{eqnarray}
&&\Pi_3(q)=-\int_{\mathbb{T}^d}\int_{\mathbb{R}}\mu(\xi,q)\int_{\mathbb{R}^d}J(z,r)\mu(\xi-z,q-r)rdzdr d\xi-\hat{\alpha}, \label{Equ3.23}\\
&&\Pi_4(\xi,q)=-\int_{\mathbb{R}}\mu(\xi,q)\int_{\mathbb{R}^d}J(z,r)\mu(\xi-z,q-r)rdzdr-\hat{\alpha}-\Pi_3, \label{Equ3.24}
\end{eqnarray}
thus we have
\begin{eqnarray}\label{Equ3.25}
\mu(\xi,q)\mathfrak{L}u_{2}(x,t,\xi,q)=(\Pi_1+\Pi_2)\cdot\nabla\nabla u+(\Pi_3+\Pi_4)\frac{\partial u}{\partial t}.
\end{eqnarray}
From Lemma \ref{lem3.2}, let $\chi_{22}(\xi,q)$ be a stationary zero mean solution of the equation
\begin{eqnarray}\label{Equ3.26}
\mu(\xi,q)\mathfrak{L}\chi_{22}(\xi,q)=\Pi_2-\Pi_4\hat{\alpha}^{-1}\Theta,
\end{eqnarray}
then $\chi_{22}$ is sub-linear growth due to its stationary, that is $\varepsilon\chi_{22}\to 0$ as $\varepsilon \to 0^{+}.$

Actually the processes $\displaystyle{\int_{0}^{s}\Pi_i(r)dr}$ $(i=1,3)$ need not be stationary, so we consider the problem
\begin{eqnarray}\label{Equ3.27}
\left\{\begin{array}{l}
L^\varepsilon\zeta^{\varepsilon,1}=\Pi_1\cdot\nabla\nabla u+\Pi_3\frac{\partial u}{\partial t}=(\Pi_1-\Pi_3\hat{\alpha}^{-1}\Theta)(\frac{t}{\varepsilon^2})\nabla\nabla u, \\
\zeta^{\varepsilon,1}(x,t)=0, t\leq0, \end{array}\right.
\end{eqnarray}
where
\begin{eqnarray}\label{Equ3.28}
L^\varepsilon v(x,t)=\frac{1}{\varepsilon^{d+4}}\int\int_{\mathbb{R}^{d}\times\mathbb{R}}J(\frac{x-y}{\varepsilon},\frac{t-s}{\varepsilon^2})\nu(\frac{x}{\varepsilon},\frac{t}{\varepsilon^2}, \frac{y}{\varepsilon}, \frac{s}{\varepsilon^2})\Big(v(y,s)-v(x,t)\Big)dyds.
\end{eqnarray}
This means the representation
\begin{eqnarray}\label{Equ3.29}
u^{\varepsilon}(x,t)=u_{0}(x,t)+\varepsilon\chi(\frac{x}{\varepsilon},\frac{t}{\varepsilon^{2}})\nabla u_{0}(x,t)-\zeta^{\varepsilon,1}-\varepsilon^{2}\chi_{22}(\frac{x}{\varepsilon},\frac{t}{\varepsilon^{2}})\nabla\nabla u_{0}(x,t)+\cdots.
\end{eqnarray}
It is straightforward to check that there exists a deterministic constant $C$ such that
\begin{eqnarray*}
\|\Pi_{i}\|_{L^{2}((0,1)\times\mathbb{T}^{d})}\leq C\;(i=2,4),
\end{eqnarray*}
therefore, according to the Lemma \ref{lem3.2}, the equation (\ref{Equ3.26}) has a stationary solution $\chi_{22}\in L^{2}([0,1], L^2(\mathbb{T}^{d})).$ By the definition, for a stationary zero mean process $\Lambda(s)$, denote
\begin{eqnarray*}
\varepsilon\chi_{21}(\frac{t}{\varepsilon^{2}})=\varepsilon\int_{0}^{\frac{t}{\varepsilon^{2}}}\Lambda(r)dr,
\end{eqnarray*}
then we have
\begin{eqnarray*}
\chi_{21}(s)=\int_{0}^{s}(\Pi_1-\Pi_3\hat{\alpha}^{-1}\Theta)^{ij}(r)dr=\int_{0}^{s}\Lambda^{ij}(r)dr.
\end{eqnarray*}
(\ref{Equ3.4}) implies that $\|\Pi_{i}\|_{L^{2}(0,1)}\leq C\;(i=1,3)$ with a deterministic constant $C$. From Lemma \ref{lem3.1} it holds that
\begin{eqnarray*}
\int_{0}^{\infty}\|\mathbb{E}\{\Lambda(a)\mid\mathcal{F}_{\leq 0}^{\Lambda}\}\|_{(L^{2}(\Omega))^{d^{2}}}da\leq C\int_{0}^{\infty}(e^{- \frac{ \varrho a}{2}}+\rho _\Lambda(\frac{a}{2}))da<\infty.
\end{eqnarray*}
Applying the invariance principle for this process a.s., in the space $(C[0,T])^{d^{2}}$ (see \cite[Theorem VIII.3.79]{Jacod1991}),
\begin{eqnarray}
&&\varepsilon\int_{0}^{\frac{t}{\varepsilon^{2}}}\Lambda(r)dr\stackrel{\varepsilon \rightarrow 0}{\longrightarrow}\Upsilon^{1/2}W,\; \forall\; T>0, \label{Equ3.30}\\
&&\Upsilon^{nmij}=\int_{0}^{\infty}\mathbb{E}\left(\Lambda^{nm}(0)\Lambda^{ij}(r)+\Lambda^{ij}(0)\Lambda^{nm}(r)\right)dr, \label{Equ3.31}
\end{eqnarray}
here $W$ is a standard $d^{2}-$dimensional Wiener process. $\Lambda$ is symmetric and negative semi-definite $d^{2}\times d^{2}$ matrix, its square root is well defined.

For smooth deterministic function $u_{0}(x,t)$, we introduce an auxiliary function $\hat{\zeta}^{\varepsilon}$ which is the solution to the following problem
\begin{eqnarray}\label{Equ3.32}
\left\{\begin{array}{l}
\varepsilon\hat{\alpha}\frac{\partial \hat{\zeta}^{\varepsilon}}{\partial t}+\varepsilon\Theta\cdot\nabla\nabla\hat{\zeta}^{\varepsilon}
=\Lambda^{ij}(\frac{t}{\varepsilon^{2}})\frac{\partial^{2}}{\partial x^{i}\partial x^{j}}u_{0}(x,t),\\
\hat{\zeta}^{\varepsilon}(x,0)=0. \end{array}\right.
\end{eqnarray}
Since $v^{ij}_{0}(x,t)=\frac{\partial^{2}}{\partial x^{i}\partial x^{j}}u_{0}(x,t)$ solves the equation $\hat{\alpha}\partial_{t}v^{ij}_{0}+\Theta\cdot\nabla\nabla v^{ij}_{0}=0,$ $i,j=1,\cdots,d,$ respectively, then we have
\begin{eqnarray}\label{Equ3.33}
\hat{\zeta}^{\varepsilon}(x,t)=\varepsilon\hat{\alpha}^{-1}\chi_{21}(\frac{t}{\varepsilon^{2}})\nabla\nabla u_0(x,t).
\end{eqnarray}
Denote
\begin{eqnarray}\label{Equ3.34}
\varepsilon^{-1}\zeta^{\varepsilon,1}=\varepsilon\hat{\alpha}^{-1}\chi_{21}(\frac{t}{\varepsilon^{2}})\nabla\nabla u_0(x,t)+\Sigma^{\varepsilon}(x,t),
\end{eqnarray}
then we have the convergence of $\Sigma^{\varepsilon}(x,t)$ as $\varepsilon\rightarrow 0^{+}$.
\begin{lemma}\label{lem3.4} As $\varepsilon\rightarrow 0^{+}$, $\Sigma^{\varepsilon}$ converges to zero in probability in $L^{2}(\mathbb{R}^{d}\times(0,T)).$
\end{lemma}
\noindent{\bf Proof.} The conclusion is obtained from $||\Sigma^{\varepsilon}||\leq C\varepsilon^\gamma, \gamma>0$, the detailed proof will be given in the last Section 5. Furthermore the condition of $J$ in Lemma \ref{lem3.4} can be weakened. $\Box$
\begin{lemma}\label{lem3.5} As $\varepsilon\rightarrow 0^{+}$, the functions $\varepsilon^{-1}\zeta^{\varepsilon,1}$ converge a.s. in the space $L^{2}((0,T),L^{2}(\mathbb{R}^{d}))$ to a unique solution of the following SPDE with a finite dimensional additive noise
\begin{eqnarray}\label{Equ3.35}
\left\{\begin{array}{l}
\hat{\alpha}d\zeta^{0,1}+\Theta\cdot\nabla\nabla\zeta^{0,1}dt=(\Upsilon^{1/2})^{ijkl}\frac{\partial^{2}}{\partial x^{i}\partial x^{j}}u_{0}(x,t)dW_{t}^{kl}, \\
\zeta^{0,1}(x,0)=0. \end{array}\right.
\end{eqnarray}
\end{lemma}
\noindent{\bf Proof.} Together (\ref{Equ3.34}) with the Lemma \ref{lem2.1}, (\ref{Equ3.28}) and (\ref{Equ3.32}) yields the following equation for $\Sigma^{\varepsilon}$:
\begin{eqnarray}\label{Equ3.36}
\left\{\begin{array}{l}
\hat{\alpha}^\varepsilon\frac{\partial\Sigma^{\varepsilon}}{\partial t}+\Theta^\varepsilon\nabla\nabla\Sigma^{\varepsilon}+\Phi_0^*(\Sigma^{\varepsilon})
=\varepsilon\hat{\alpha}^{-1}\chi_{21}(\frac{t}{\varepsilon^{2}})(\Theta-\Theta^\varepsilon)\nabla\nabla v_{0}\\
\qquad\quad+\varepsilon\hat{\alpha}^{-1}\partial_{t}\chi_{21}(\frac{t}{\varepsilon^{2}})(\hat{\alpha}-\hat{\alpha}^\varepsilon)v_{0}+\varepsilon\hat{\alpha}^{-1}\chi_{21}(\frac{t}{\varepsilon^{2}}) (\hat{\alpha}-\hat{\alpha}^\varepsilon)\partial_{t} v_{0}-\Phi_0^*(\hat{\zeta}^{\varepsilon}),\\
\Sigma^{\varepsilon}(x,0)=0, \end{array}\right.
\end{eqnarray}
where
\begin{eqnarray}
&&\hat{\alpha}^\varepsilon=\int_{\mathbb{R}}\mu(\frac{x}{\varepsilon},\frac{t}{\varepsilon^2})\int_{\mathbb{R}^d}J(z,r)\mu(\frac{x}{\varepsilon}-z,\frac{t}{\varepsilon^2}-r)(-r)dzdr, \label{Equ3.37}\\
&&\Theta^\varepsilon={\frac{1}{2}}\int_{\mathbb{R}}\mu(\frac{x}{\varepsilon},\frac{t}{\varepsilon^2})\int_{\mathbb{R}^d}J(z,r)\mu(\frac{x}{\varepsilon}-z,\frac{t}{\varepsilon^2}-r)z\otimes zdzdr,\label{Equ3.38}\\
&&\Phi_0^{*}(u)(x,t)=\int_{\mathbb{R}}\mu(\frac{x}{\varepsilon},\frac{t}{\varepsilon^2})\int_{\mathbb{R}^{d}}J(z,r)\mu(\frac{x}{\varepsilon}-z,\frac{t}{\varepsilon^2}-r) \nonumber\\
&&\cdot \Big\{\Big[\int_{0}^{1}\nabla\nabla u(x-\varepsilon z\theta,t-\varepsilon^2 r\theta)\cdot z\otimes z(1-\theta)d\theta-\frac{1}{2}\nabla \nabla u(x,t)\cdot z\otimes z\Big]\nonumber\\
&&+2\varepsilon r\int_{0}^{1}\partial_t\nabla u(x-\varepsilon z\theta,t-\varepsilon^2 r\theta)\otimes z(1-\theta)d\theta\nonumber\\
&&+\varepsilon^2 r^2\int_{0}^{1}\partial_{tt} u(x-\varepsilon z\theta,t-\varepsilon^2 r\theta)(1-\theta)d\theta\Big\}dzdr. \label{Equ3.39}
\end{eqnarray}
When $\varepsilon \to 0,$ from regularity of $u_0$ and (\ref{Equ3.33}), (\ref{Equ3.39}), it's easy to deduce that $\Phi_0^*(\hat{\zeta}^{\varepsilon})\rightharpoonup0$ a.s. in $L^2((0,T)\times \mathbb{R}^d)$, actually with respect to the spatial variable it weakly converges to its average in $L^2((0,T)\times Q)$ and with respect to the time variable it weakly converges to an expectation in $L^2(0,T)$ by the Birkhoff ergodic theorem, then this also yields in the sense of $L^2((0,T)\times L_{loc}^2 (\mathbb{R}^d))$
\begin{eqnarray}
&&\Theta^\varepsilon\rightharpoonup\frac{1}{2}\mathbb{E}\int_{\mathbb{T}^{d}}\int_{\mathbb{R}}\mu(\xi,q)\int_{\mathbb{R}^{d}}J(z,r)\mu(\xi-z,q-r)z^2dzdrd\xi, \label{Equ3.40}\\
&&\hat{\alpha}^\varepsilon\rightharpoonup-\mathbb{E}\int_{\mathbb{T}^{d}}\int_{\mathbb{R}}\mu(\xi,q)\int_{\mathbb{R}^{d}}J(z,r)\mu(\xi-z,q-r)rdzdrd\xi=\hat{\alpha}. \label{Equ3.41}
\end{eqnarray}
Since $\hat{\alpha}-\hat{\alpha}^\varepsilon,\Theta-\Theta^\varepsilon$ are bounded in $L^2((0,T)\times \mathbb{R}^d)$ and $\varepsilon\chi_{2,1}(\frac{t}{\varepsilon^{2}})$ is a compact set in $C(0,T)$, then there exists a subsequence that $\Sigma^{\varepsilon}$ converges to some function $\Sigma_0$ from \cite[Lemmas 5.3-5.5]{PiatnitskiStochastic}. Noticed that $\varepsilon\chi_{2,1}(\frac{t}{\varepsilon^{2}})$ need not tend to zero in $[0,T]$, next we will prove  $\Sigma^{\varepsilon}$ tends to zero in $L^2(\mathbb{R}^d \times(0,T))$ as $\varepsilon \to 0.$ Due to the difference in the form of the equations (\ref{Equ3.27}) and (\ref{Equ3.32}), we can't be directly investigate the limiting behaviour of  $\Sigma^{\varepsilon}$ in (\ref{Equ3.36}), the difficulty of the problem lies in finding mesoscopic states between nonlocal operators and local operators.

Passing the limit in (\ref{Equ3.33}) we obtain that $\hat{\zeta}^{\varepsilon}$ converges in law in $C(0,T; L^{2}(\mathbb{R}^{d}))$ to the process $\Upsilon^{1/2}W$ from (\ref{Equ3.30}) and (\ref{Equ3.31}). Lemma \ref{lem3.4} shows that the limits of $\varepsilon^{-1}\zeta^{\varepsilon,1}$ and $\hat{\zeta}^{\varepsilon}$ in $L^2(\mathbb{R}^d\times(0,T))$ are equal. The convergence result and Lemma \ref{lem3.5} is proved. $\Box$

We also require that $u^{\varepsilon}$ satisfies the initial condition at the level $(\varepsilon^{1})$. In order to ensure that initial value of expansion function will not be influenced by $\chi$, we introduce one more term of order $\varepsilon^{1}$ so that the expansion takes the form $\mathcal{I}$ satisfying that
\begin{eqnarray}\label{Equ3.42}
\left\{\begin{array}{l}
\displaystyle{\int_{\mathbb{R}}\mu(x,t)\int_{\mathbb{R}^d}J(z,r)\mu(x-z,t-r)[\mathcal{I}(x-z,t-r)-\mathcal{I}(x,t)]dzdr=0,}\\
\mathcal{I}(x,t)=-\chi(x,0), t\leq0, \end{array}\right.
\end{eqnarray}
for $(x,t)\in \mathbb{T}^d\times\mathbb{R}^+$, we can deduce that the solution $\mathcal{I}(x,t)$ decays exponentially as $t\rightarrow +\infty$ just like $\chi$ and for all $1<p<+\infty$ and positive constants $C, \gamma$ such that
\begin{eqnarray}\label{Equ3.43}
\int_{\mathbb{T}^{d}}\mathcal{I}(z,t)dz=\int_{\mathbb{T}^{d}}\mathcal{I}(z,0)dz=0,\;\|\mathcal{I}(\cdot,t)\|_{L^{p}(\mathbb{T}^{d})}\leq Ce^{-\gamma t}, t>0.
\end{eqnarray}
The final multi-scale expansion is as
\begin{eqnarray}\label{Equ3.44}
u^{\varepsilon}(x,t)&=&u_{0}(x,t)+\varepsilon\Big(\chi(\frac{x}{\varepsilon},\frac{t}{\varepsilon^{2}})+\mathcal{I}(\frac{x}{\varepsilon},\frac{t}{\varepsilon^{2}})\Big)\nabla u_0(x,t)\nonumber\\
&&-\zeta_1^\varepsilon-\varepsilon^{2}\chi_{22}(\frac{x}{\varepsilon},\frac{t}{\varepsilon^{2}})\nabla\nabla u_0+\cdots.
\end{eqnarray}

Finally we consider terms in $(\varepsilon^{1})$. Its right-hand side can be rearranged as follows:
\begin{eqnarray}\label{Equ3.45}
&&\int_{\mathbb{R}}\mu(\frac{x}{\varepsilon},\frac{t}{\varepsilon^2})\int_{\mathbb{R}^{d}}J(z,r)\mu(\frac{x}{\varepsilon}-z,\frac{t}{\varepsilon^2}-r)\nonumber\\
&&\cdot\Big\{\int_{0}^{1}\nabla\nabla u(x-\varepsilon z\theta,t-\varepsilon^2 r\theta)\cdot z\otimes z(1-\theta)d\theta-\frac{1}{2}\nabla\nabla u\cdot z\otimes z \nonumber\\
&&+\varepsilon\chi(\frac{x}{\varepsilon}-z,\frac{t}{\varepsilon^2}-r)\int_{0}^{1}\nabla\nabla\nabla u(x-\varepsilon z\theta,t-\varepsilon^2 r\theta)z\otimes z(1-\theta)d\theta \nonumber\\
&&-\varepsilon\chi_{22}(\frac{x}{\varepsilon}-z,\frac{t}{\varepsilon^2}-r)\int_{0}^{1}\nabla\nabla\nabla u(x-\varepsilon z\theta,t-\varepsilon^2 r\theta)zd\theta\Big\}dzdr\nonumber\\
&&-\varepsilon\int_{\mathbb{R}}\mu(\frac{x}{\varepsilon},\frac{t}{\varepsilon^2})\int_{\mathbb{R}^d}J(z,r)\mu(\frac{x}{\varepsilon}-z,\frac{t}{\varepsilon^2}-r)
r\chi(\frac{x}{\varepsilon}-z,\frac{t}{\varepsilon^2}-r)\frac{\partial\nabla u}{\partial t}dzdr\nonumber\\
&&+2\varepsilon \int_{\mathbb{R}}\mu(\frac{x}{\varepsilon},\frac{t}{\varepsilon^2})\int_{\mathbb{R}^d}J(z,r)\mu(\frac{x}{\varepsilon}-z,\frac{t}{\varepsilon^2}-r)r\nonumber\\
&&\cdot\int_{0}^{1}\partial_t\nabla u(x-\varepsilon z\theta,t-\varepsilon^2 r\theta)\otimes z(1-\theta)d\theta dzdr:=\varepsilon\Psi_{3}(\frac{x}{\varepsilon},\frac{t}{\varepsilon^{2}},u).
\end{eqnarray}
Using Taylor's expansion, it is easy to check that
\begin{eqnarray}\label{Equ3.46}
&&\Psi_{3}(\frac{x}{\varepsilon},\frac{t}{\varepsilon^{2}},u)=\int_{\mathbb{R}}\mu(\frac{x}{\varepsilon},\frac{t}{\varepsilon^{2}})\int_{\mathbb{R}^{d}}J(z,r)\mu(\frac{x}{\varepsilon}-z,\frac{t}{\varepsilon^{2}}-r)\nonumber\\
&&\qquad\cdot\Big[\Big(-\frac{1}{6}z^3+\frac{1}{2}z^2 \chi(\frac{x}{\varepsilon}-z,\frac{t}{\varepsilon^{2}}-r)-\chi_{22}(\frac{x}{\varepsilon}-z,\frac{t}{\varepsilon^{2}}-r)\Big)\frac{\partial^{3}u}{\partial x^{3}}\nonumber\\
&&\qquad+r(z-\chi(\frac{x}{\varepsilon}-z,\frac{t}{\varepsilon^2}-r))\frac{\partial\nabla u}{\partial t}\Big]dzdr+o(\varepsilon)\nonumber\\
&&\qquad:=\Psi_{4}(\frac{x}{\varepsilon},\frac{t}{\varepsilon^{2}})\frac{\partial^{3}u}{\partial x^{3}}+\Psi_{5}(\frac{x}{\varepsilon},\frac{t}{\varepsilon^{2}})\frac{\partial\nabla u}{\partial t}+o(\varepsilon),
\end{eqnarray}
where $\frac{\partial^{3}u}{\partial x^{3}}(x,t)=(\frac{\partial^{3}u}{\partial x^{i}\partial x^{j}\partial x^{k}})_{i,j,k=1}^{d}$ stands for the tensor of third order partial derivatives of $u$, $o(\varepsilon)$ represent higher-order remainders of Taylor expansions.

We first introduce the following constant tensor $\mu=\{\mu^{ijk}\}_{i,j,k=1}^{d}$ and constant vector $\mu_1$:
\begin{eqnarray*}
\mu&=&\mathbb{E}\int_{\mathbb{T}^{d}}\int_{\mathbb{R}}\mu(\xi,q)\int_{\mathbb{R}^{d}}J(z,r)\mu(\xi-z,q-r)\\
&&\cdot(-\frac{1}{6}z^3+\frac{1}{2}z^2\chi(\xi-z,q-r)-\chi_{22}(\xi-z,q-r))dzdrd\xi, \\
\mu_1&=&\mathbb{E}\int_{\mathbb{T}^{d}}\int_{\mathbb{R}}\mu(\xi,q)\int_{\mathbb{R}^{d}}J(z,r)\mu(\xi-z,q-r)r(z-\chi(\xi-z,q-r))dzdrd\xi.
\end{eqnarray*}
Then we consider the following problems
\begin{eqnarray}
&&\left\{\begin{array}{l}
L^\varepsilon\vartheta_{\varepsilon,1}=\Psi_{3}(\frac{x}{\varepsilon},\frac{t}{\varepsilon^{2}},u)-\mu\frac{\partial^{3}u}{\partial x^{3}}-\mu_1\frac{\partial\nabla u}{\partial t}, \\
\vartheta_{\varepsilon,1}(x,t)=0, t\leq0, \end{array}\right. \label{Equ3.47}\\
&&\left\{\begin{array}{l}
L^\varepsilon\vartheta_{\varepsilon,2}=\mu\frac{\partial^{3}}{\partial x^{3}}u(x,t)+\mu_1\frac{\partial\nabla u}{\partial t}, \\
\vartheta_{\varepsilon,2}(x,t)=0, t\leq0. \end{array}\right. \label{Equ3.48}
\end{eqnarray}

First we give an important lemma.
\begin{lemma}\label{lem3.6} $\hbar^\varepsilon(x,t)=\hbar(\frac{x}{\varepsilon},\frac{t}{\varepsilon^2})\in L^2_{loc}(\mathbb{R}, L^2_{per}(\mathbb{T}^d))$ is periodic in $x$ and stationary in $t$, for any cube
$U\subset\mathbb{R}^{d}$, we have
\begin{eqnarray}\label{Equ3.49}
\hbar(\frac{x}{\varepsilon},\frac{t}{\varepsilon^2})\rightharpoonup\mathbf{E}\int_{\mathbb{T}^{d}}\hbar(z,s)dz=\hat{\mu}\;a.s. ~\mbox{as}\;\varepsilon\rightarrow 0^{+}~weakly ~in ~ L^{2}(U\times(0,T)).
\end{eqnarray}
\end{lemma}
\noindent{\bf Proof.} We decompose the function $\hbar(z,s)$ as
\begin{eqnarray*}
\hbar(z,s)=\mathcal{M}_Y{\hbar}(s)+\left[\hbar(z,s)-\mathcal{M}_Y{\hbar}(s)\right],\;\mathcal{M}_Y{\hbar}(s)=\int_{\mathbb{T}^{d}}\hbar(z,s)dz.
\end{eqnarray*}
Then for all $s, \hbar-\mathcal{M}_Y{\hbar} \in {L^{2}(\mathbb{T}^{d}\times(s,s+1))}$ and $\displaystyle{\int_{\mathbb{T}^{d}}(\hbar(x,s)-\mathcal{M}_Y{\hbar}(s))dx=0.}$ Therefore, $\hbar^{\varepsilon}-(\mathcal{M}_Y{\hbar})^{\varepsilon}$ converges weakly to zero in $L^{2}(U\times(0,T))$ as $\varepsilon\rightarrow 0$. From $\mathcal{M}_Y{\hbar}$ and $\hat{\mu}$ it can be seen directly  that $\mathcal{M}_Y{\hbar}$ is stationary, and $\mathbb{E}\left(\mathcal{M}_Y{\hbar}(t)-\hat{\mu}\right)=0$. $\left((\mathcal{M}_Y{\hbar})^{\varepsilon}-\hat{\mu}\right)$ converges a.s. to zero weakly in $L^{2}(0,T)$ by the Birkhoff ergodic theorem. The function $\left(\hbar^{\varepsilon}-\hat{\mu}\right)$ converges a.s. to zero weakly in $L^{2}(U\times$ $(0,T))$. This completes the proof of Lemma \ref{lem3.6}. $\Box$

Due to the regularity of $u_0$ and periodicity of $\hbar$ in spatial variable this implies that $$\left(\hbar\left(\frac{x}{\varepsilon}, \frac{t}{\varepsilon^{2}}\right)-\mu\right)\frac{\partial^{3}}{\partial x^{3}}u_{0}$$ converges a.s. to zero weakly in $L^{2}\left(\mathbb{R}^{d}\times(0,T)\right)$.
\begin{lemma}\label{lem3.7} The solution of problem (\ref{Equ3.47}) tends to zero a.s. as $\varepsilon\rightarrow 0^{+}$ in $L^{2}\left(Q\times[0,T]\right)$ for any cube $Q$. Moreover,
$\lim\limits_{\varepsilon\rightarrow 0}\mathbb{E}\left(\left\|\vartheta_{\varepsilon, 1}\right\|_{L^{2}\left(Q\times[0, T]\right)}^{2}\right)=0.$
\end{lemma}
\noindent{\bf Proof.} We first prove that sequence $\{\vartheta_{\varepsilon,1}\}_{\varepsilon>0}$ is bounded, then it possess compactness, finally we prove that the sequence converges to zero. Since
\begin{eqnarray}
&&\Psi_{4}(\frac{x}{\varepsilon},\frac{t}{\varepsilon^{2}})=\int_{\mathbb{R}}\mu(\frac{x}{\varepsilon},\frac{t}{\varepsilon^{2}})\int_{\mathbb{R}^{d}}J(z,r)  \mu(\frac{x}{\varepsilon}-z,\frac{t}{\varepsilon^{2}}-r)\nonumber\\
&&\qquad\cdot\Big(-\frac{z^3}{6}+\frac{z^2}{2}\chi(\frac{x}{\varepsilon}-z,\frac{t}{\varepsilon^{2}}-r)-\chi_{22}(\frac{x}{\varepsilon}-z,\frac{t}{\varepsilon^{2}}-r)\Big)dzdr, \label{Equ3.50}\\
&&\Psi_{5}(\frac{x}{\varepsilon}, \frac{t}{\varepsilon^{2}})-\mu_1=\int_{\mathbb{R}}\mu(\frac{x}{\varepsilon},\frac{t}{\varepsilon^2}) \int_{\mathbb{R}^d}J(z,r)\mu(\frac{x}{\varepsilon}-z,\frac{t}{\varepsilon^2}-r)r\nonumber\\
&&\qquad\cdot\Big(z-\chi(\frac{x}{\varepsilon}-z,\frac{t}{\varepsilon^2}-r)\Big)dzdr-\mu_1, \label{Equ3.51}
\end{eqnarray}
from integrability of $J$ we conclude that $\|\Psi_{i}^{\varepsilon}\|_{L^{2}([0,T]\times\mathbb{R}^{d})}\leq C\varepsilon\;(i=4,5)$, where $\Psi^{\varepsilon}_{i}=\Psi_{i}(\frac{x}{\varepsilon},\frac{t}{\varepsilon^{2}})$, $J^\varepsilon(x,t)=J(\frac{x}{\varepsilon},\frac{t}{\varepsilon^{2}})$ and $\chi_{22}^{\varepsilon}(x, t)=\chi_{22}(\frac{x}{\varepsilon},\frac{t}{\varepsilon^{2}})$. Due to the regularity of $u_{0}$, with a deterministic $C$, we have
\begin{eqnarray*}
\Big\|\Psi_{3}(\frac{x}{\varepsilon},\frac{t}{\varepsilon^{2}},u)-\mu\frac{\partial^{3}u}{\partial x^{3}}-\mu_1\frac{\partial\nabla u}{\partial t}\Big\|_{L^{2}(\mathbb{R}^{d}\times(0,T))}\leq C.
\end{eqnarray*}

We will show that the family $\vartheta_{\varepsilon, 1}$ is compact in $L^{2}((0,T),L_{{loc}}^{2}(\mathbb{R}^{d})),$ first we decompose $\vartheta_{\varepsilon,1}=\vartheta^1_{\varepsilon,1} +\vartheta^2_{\varepsilon,1}$, where $\vartheta^1_{\varepsilon,1}$ and $\vartheta^2_{\varepsilon,1}$ are solutions of the following problems respectively
\begin{eqnarray}
&&\left\{\begin{array}{l}
L^\varepsilon\vartheta^1_{\varepsilon,1}=\Psi_{4}\frac{\partial^{3}}{\partial x^{3}}u(x,t)+\Psi_{5}\frac{\partial\nabla u}{\partial t}-\mu\frac{\partial^{3}u}{\partial x^{3}}-\mu_1\frac{\partial\nabla u}{\partial t}, \\
\vartheta^1_{\varepsilon,1}(x,t)=0, t\leq0, \end{array}\right. \label{Equ3.52}\\
&&\left\{\begin{array}{l}
L^\varepsilon\vartheta^2_{\varepsilon,1}=o(\varepsilon), \\
\vartheta^2_{\varepsilon,1}(x,t)=0,t\leq0. \end{array}\right. \label{Equ3.53}
\end{eqnarray}
Piatnitski and Zhizhina \cite[Lemmas 5.3-5.5]{PiatnitskiStochastic} proved the analogous compactness by constructing the finite 2$\delta$ net for $\delta>0$, they provided a very skillful proof by the fact that $\vartheta^1_{\varepsilon,1}$ is compact and  $\vartheta^2_{\varepsilon,1}$ converges to zero as $\varepsilon \to 0$ implies the compactness of the family $\{\vartheta_{\varepsilon,1}\}$.

Due to the Lemma \ref{lem3.6}, the function $\Psi_{3}^{\varepsilon}-\mu$ converges a.s. to zero weakly in $L^{2}(Q\times$ $(0,T))$ for any cube $Q\subset\mathbb{R}^{d}$. By the Lebesgue dominated convergence theorem and combining this with the above compactness arguments, we conclude that a.s. $\vartheta_{\varepsilon,1}$ converges to zero in $L^{2}\left(Q\times(0,T)\right)$. This completes the proof of Lemma \ref{lem3.7}. $\Box$

We now present the main result in this section. Denote
\begin{eqnarray}\label{Equ3.54}
U^{\varepsilon}(x,t)=u^{\varepsilon}(x,t)-u_{0}(x,t)-\varepsilon\chi\Big(\frac{x}{\varepsilon},\frac{t}{\varepsilon^{2}}\Big)\nabla_{x}u_{0},\; \omega^\varepsilon=\varepsilon^{-1}U^{\varepsilon}.
\end{eqnarray}
The limit behaviour of $\omega^{\varepsilon}$ is given by the following theorem.
\begin{theorem}\label{the3.2} The function $\omega^{\varepsilon}$ converges in law as $\varepsilon\rightarrow 0^{+}$ in $L^{2}\left(\mathbb{R}^{d}\times(0,T)\right)$ to a unique solution of the following SPDE
\begin{eqnarray}\label{Equ3.55}
\hat{\alpha}d\omega_{0}+\Theta\cdot\nabla\nabla\omega_{0}dt=-(\mu\frac{\partial^{3}u}{\partial x^{3}}+\mu_1\partial_t\partial_xu)dt-\Upsilon^{1/2}\frac{\partial^{2}u}{\partial x^{2}}dW_{t}.
\end{eqnarray}
\end{theorem}
\noindent{\bf Proof.} We set
\begin{eqnarray}\label{Equ3.56}
\mathcal{B}^{\varepsilon}(x,t)&=&\omega^{\varepsilon}+\varepsilon^{-1}\zeta^{\varepsilon,1}(x,t)+\vartheta_{\varepsilon,2}(x,t)-\mathcal{I}(\frac{x}{\varepsilon},\frac{t}{\varepsilon^{2}})\nabla u(x,t)\nonumber \\
&+&\varepsilon\chi_{22}(\frac{x}{\varepsilon},\frac{t}{\varepsilon^{2}})\nabla\nabla u(x,t)+\vartheta_{\varepsilon,1}(x,t).
\end{eqnarray}
Substituting $\mathcal{B}^{\varepsilon}$ in (\ref{Equ3.1}) for $u^{\varepsilon}$ and combining the above equations, after simple calculation and rearrangement according to the order of $\varepsilon$ we obtain that $\mathcal{B}^{\varepsilon}$ satisfies the problem
\begin{eqnarray}\label{Equ3.57}
L^\varepsilon \mathcal{B}^{\varepsilon}=R^{\varepsilon},\; 1_{\{t\leq0\}}\mathcal{B}^{\varepsilon}(x,t)=R_{1}^{\varepsilon},
\end{eqnarray}
where
\begin{eqnarray}\label{Equ3.58}
R^{\varepsilon}&=&-\int_{\mathbb{R}}\mu(\frac{x}{\varepsilon},\frac{t}{\varepsilon^2})\int_{\mathbb{R}^d}J(z,r)\mu(\frac{x}{\varepsilon}-z,\frac{t}{\varepsilon^2}-r)\nonumber\\
&&\cdot\Big\{\varepsilon^{-2}[\mathcal{I}(\frac{x}{\varepsilon}-z,\frac{t}{\varepsilon^2}-r)-\mathcal{I}(\frac{x}{\varepsilon},\frac{t}{\varepsilon^2})]\partial_{x}u-\varepsilon^{-1}[z\otimes \mathcal{I}(\frac{x}{\varepsilon}-z,\frac{t}{\varepsilon^2}-r)]\partial_{x}^{2}u\nonumber\\
&&-r\mathcal{I}(\frac{x}{\varepsilon}-z,\frac{t}{\varepsilon^2}-r)\partial_{t}\partial_{x}u+[z\otimes\chi_{22}(\frac{x}{\varepsilon}-z,\frac{t}{\varepsilon^2}-r)]\partial_{x}^{3}u\nonumber\\
&&+\varepsilon[z^2\otimes\chi_{22}(\frac{x}{\varepsilon}-z,\frac{t}{\varepsilon^2}-r)]\partial_{x}^{4}u
+r\varepsilon\chi_{22}(\frac{x}{\varepsilon}-z,\frac{t}{\varepsilon^2}-r)\partial_{t}\partial_{x}^{2}u\Big\}dzdr\nonumber\\
&&+o^{\mathcal{I}}(\varepsilon)+o^{\chi_{22}}(\varepsilon):=I_1+I_2+\cdots+I_8,
\end{eqnarray}
here $I_7, I_8 $ are the higher order remainders of using Taylor expansion formulas. Due to (\ref{Equ3.42}), $I_1=0$. Since $\|\mathcal{I}(\cdot,t)\|_{L^{p}(\mathbb{T}^{d})}\leq Ce^{-\gamma t}$, $I_3$ tends to zero as $\varepsilon \to 0.$
It's not hard to verify that $|I_i|\leq C\varepsilon$ except $i=2,4$. Let $\vartheta_{\varepsilon,3}$ be the solution of
\begin{eqnarray}\label{Equ3.59}
L^\varepsilon\vartheta_{\varepsilon,3}=I_2, \quad \vartheta_{\varepsilon,3}(x,t)=0,t\leq0.
\end{eqnarray}
Next we prove that $\vartheta_{\varepsilon,3}$ is compact (similar process with $\vartheta_{\varepsilon,1}$) and $I_2$ converges a.s. to zero in $L^{2}((0,T)\times\mathbb{R}^{d})$ and
\begin{eqnarray}\label{Equ3.60}
\lim_{\varepsilon\rightarrow 0}\mathbb{E}(\|\vartheta_{\varepsilon,3}\|_{L^{2}(\mathbb{R}^{d}\times(0,T))}^{2})=0.
\end{eqnarray}
Similarly we can get the estimate of $I_4$. We have known $\varepsilon \chi_{22}\left(\frac{x}{\varepsilon},\frac{t}{\varepsilon^{2}}\right),\vartheta_{\varepsilon,1}(x, t), \mathcal{I}\left(\frac{x}{\varepsilon}, \frac{t}{\varepsilon^{2}}\right)\to 0$ as $\varepsilon \to 0$ in $ L^{2}\left(Q\times(0,T)\right)$ from (\ref{Equ3.26}), Lemma \ref{lem3.6} and the definition of $\mathcal{I}$. Combining the above estimates and the  regularity of $u$ from condition $\mathbb{P}$ we conclude that $R^{\varepsilon}$ a.s. tends to zero in $L^{2}(\mathbb{R}^{d}\times(0,T))$, as $\varepsilon\rightarrow 0$, and $R_{1}^{\varepsilon}$ a.s. tends to zero in $L^{2}((-\infty,0),L^{2}(\mathbb{R}^{d}))$.

By Lemma \ref{lem3.3} and Lemma \ref{lem3.5}, it follows that $\omega^{\varepsilon}+\varepsilon^{-1}\zeta^{\varepsilon,1}+\vartheta_{\varepsilon,2}$ satisfies the convergence
$\omega^{\varepsilon}\rightarrow -\zeta^{0,1}-\vartheta_{0,2}\;\mbox{in}\;L^{2}(\mathbb{R}^{d}\times[0,T]),$ where $\vartheta_{0,2}$ satisfies the following equation
\begin{eqnarray}\label{Equ3.61}
\left\{\begin{array}{l}
\hat{\alpha}\frac{\partial\vartheta_{0,2}}{\partial t}+\Theta\cdot\nabla\nabla\vartheta_{0,2}=\mu\frac{\partial^{3}}{\partial x^{3}}u(x,t)+\mu_1\frac{\partial\nabla u}{\partial t}, \\
\vartheta_{0,2}(x,0)=0. \end{array}\right.
\end{eqnarray}
Due to $\omega_0=-\zeta^{0,1}-\vartheta_{0,2},$ this completes the proof of Theorem \ref{the3.2}. $\Box$
\begin{remark}\label{rem3.3} If $J(x,t)$ satisfies (\ref{Equ1.3}) and $\int_{\mathbb{R}^d}|x|^3J(x,t)dx<\infty, t\in \mathbb{R},$ $u^\varepsilon(x,t)$ and $u(x,t)$ are the solutions of (\ref{Equ3.1}) and (\ref{Equ3.13}) then we have the following estimate
\begin{eqnarray}\label{Equ3.62}
||u^\varepsilon-u||_{L^{2}(\mathbb{R}^{d}\times(0,T))}\leq C\varepsilon.
\end{eqnarray}
Recently \cite{Piatnitski2022} discussed the problem of nonlocal operator convergence rate with the help of the spectral edge of the operator when $J(x)$ satisfies $\int_{\mathbb{R}^d}J(z)|z|^kdz<\infty, k=1,2,3$. The above results can also support their conclusion from another aspect when $J$ compactly supported in the set $D_{\delta\gamma}$.
\end{remark}

\section{$\mu(x,t)$ is stationary in $x,t$}\label{sec4}
\setcounter{equation}{0}

We denote by $\mathbf{L}^{2}$ the set of stationary maps $u=u(x,t,\omega)$ such that
\begin{eqnarray*}
u(x+k,t+s,\omega)=u(x,t,T_{(k,s)}\omega), \forall ~(k,s,\omega)\in\mathbb{R}^{d}\times\mathbb{R}\times\Omega.
\end{eqnarray*}
Define $\mathbf{L}^{2}-$norm: $\|u\|_{\mathbf{L}^{2}}=\mathbb{E}[\int_{(-1/2,1/2)^{d+1}}u^{2}]<+\infty.$ Note that, if $u\in\mathbf{L}^{2}$ and $U$ is a bounded measurable subset of $\mathbb{R}^{d}$, the stationarity in time implies that the limit
\begin{eqnarray*}
\mathbb{E}\int_{U}u(x,t)dx=\lim_{h\rightarrow 0^{+}}\mathbb{E}\Big[\frac{1}{2h}\int_{U}\int_{t-h}^{t}u(x,s)dxds\Big].
\end{eqnarray*}
Setting
\setcounter{equation}{0}
\begin{eqnarray}\label{Equ4.1}
\mathfrak{A}_\omega u&=&\int_{\mathbb{R}}\int_{\mathbb{R}^d}J(z,r)\nu(x,t;x-z,t-r)\Big|-z+\chi(x-z,t-r)-\chi(x,t)\Big|^{p-2}\nonumber\\
&&\cdot\Big(-z+\chi(x-z,t-r)-\chi(x,t)\Big)dzdr,
\end{eqnarray}
$\nu(x,t;x-z,t-r)=\mu(x,t,\omega)\mu(x-z,t-r,\omega),$ $\nu_1(x+k,t+s,\omega)=\nu_1(x,t,T_{k+s}\omega),$ $\forall x\in \mathbb{R}^d,$ $s\in\mathbb{R}$ and
\begin{eqnarray}\label{Equ4.2}
\zeta_{z,r}(x,t,\omega)=\chi(x+z,t+r,\omega)-\chi(x,t,\omega),
\end{eqnarray}
and there exist a constant $c>0$ and a cube $\mathbf{U}\subset\mathbb{R}^{d}\times \mathbb{R}$ such that $J(x,t)\geq c$ for all $(x,t)\in\mathbf{U}$. This additional condition on $J(x,t)$ is naturally satisfied for regular kernels.
\begin{theorem}\label{the4.1} Assume above conditions are satisfied, there exists a unique map $\chi$: $\mathbb{R}^{d+1}\times\Omega\rightarrow\mathbb{R}$ such that
\begin{eqnarray}\label{Equ4.3}
\mathfrak{A}_\omega\chi^{k}=0, k=1,2,\cdots,d,
\end{eqnarray}
and for all $z,r\in\mathbb{R}^d,$ $\zeta_{z,r}(\xi,t,\omega)$ is a stationary field for any given $z,r$ and satisfy
\begin{eqnarray*}
\int_{(-1/2,1/2)^{d+1}}\zeta_{z,r}(x,t,\omega)dxdt=0~\mathbb{P}-a.s.,\;\zeta_{z}\in\mathbf{L}^{\mathbf{2}}.
\end{eqnarray*}
Moreover, as $\varepsilon\rightarrow 0$ and $\mathbb{P}-$a.s. and in expectation, for any $\psi\in\mathbb{R}^d$
\begin{eqnarray*}
\chi^{\varepsilon}(x,t,\psi,\omega)=\varepsilon\chi^\psi(\frac{x}{\varepsilon},\frac{t}{\varepsilon^{2}},\omega)\rightarrow 0\;\mbox{in}\;L_{loc}^{2}(\mathbb{R}^{d+1}).
\end{eqnarray*}
\end{theorem}

Actually $\chi(x,t)$ and $\psi$ satisfy
\begin{eqnarray}\label{Equ4.4}
&&\int_{\mathbb{R}}\int_{\mathbb{R}^d}J(z,r)\mu(\xi-z,q-r)\Big|\Big(-z+\chi(\xi-z,q-r)-\chi(\xi,q)\Big)\cdot\psi\Big|^{p-2}\nonumber\\
&&\qquad\times\Big(-z+\chi(\xi-z,q-r)-\chi(\xi,q)\Big)\cdot\psi dzdr=0,
\end{eqnarray}
and construct auxiliary function
\begin{eqnarray*}
w^\varepsilon(x,t)=u(x,t)+\varepsilon\chi(\frac{x}{\varepsilon},\frac{t}{\varepsilon^2},\omega)\nabla u(x,t),
\end{eqnarray*}
where $\chi$ need not be a stationary field and different from (\ref{Equ3.12}), but $\chi(\xi-z,q-r)-\chi(\xi,q)$ are stationary respected to $\xi,q$ for any $(z,r)\in \mathbb{R}^d\times \mathbb{R}$ .
Substituting the expression  of $w^\varepsilon$ in (\ref{Equ2.1}), we can also collect all power-like terms in (\ref{Equ2.1}) and $\varepsilon^0$ term is
\begin{eqnarray}\label{Equ4.5}
&&(p-1)\int_{\mathbb{R}}\mu(\frac{x}{\varepsilon},\frac{t}{\varepsilon^2})\int_{\mathbb{R}^d}J(z,r)\mu(\frac{x}{\varepsilon}-z,\frac{t}{\varepsilon^2}-r)\nonumber\\
&&\qquad\times\Big|\Big(-z+\chi(\frac{x}{\varepsilon}-z,\frac{t}{\varepsilon^2}-r)-\chi(\frac{x}{\varepsilon},\frac{t}{\varepsilon^2})\Big)\cdot\nabla u(x,t)\Big|^{p-2}\nonumber\\
&&\qquad\times\Big\{\Big[-z\otimes\chi(\frac{x}{\varepsilon}-z,\frac{t}{\varepsilon^2}-r)+{\frac{1}{2}}z^{2}\Big]\cdot\nabla\nabla u-r\frac{\partial u(x,t)}{\partial t}\Big\}dzdr.
\end{eqnarray}

We give the direct result due to  the stationarity of the nonlinear term.
\begin{theorem}\label{the4.2} If $J: \mathbb{R}^{d+1}\rightarrow\mathbb{R}$ is compactly supported in $D_{\delta\gamma}$ defined in (\ref{Equ1.7}). The solutions $u^{\varepsilon}(x,t)$ to $\mathbb{S}(J_{\varepsilon},\nu,u^{\varepsilon},f)=0$ converge along subsequence uniformly  on compact sets to a viscosity solution $u(x,t)$ a.s.  as $\varepsilon\rightarrow 0^{+}$ to the problem
\begin{eqnarray}\label{Equ4.6}
\mathfrak{N}(\nabla u)\frac{\partial u}{\partial t}=\mathfrak{P}(\nabla u,\nabla \nabla u), \quad u(x,0)=f(x),
\end{eqnarray}
where $\mathfrak{N}$ and $\mathfrak{P}$ depending on $J$ that are given by
\begin{eqnarray}
\mathfrak{N}(\nabla u)&=&\mathbf{E}\int_{\mathbb{R}}\int_{\mathbb{R}^d}J(z,r)\mu(0,0)\mu(-z,-r)\nonumber\\
&&\cdot\Big|\nabla u(x,t)\Big(-z+\zeta_{-z,-r}(0,0,\omega)\Big)\Big|^{p-2}rdzdr, \label{Equ4.7}\\
\mathfrak{P}(\nabla u,\nabla \nabla u)&=&\frac{1}{2}\mathbf{E}\int_{\mathbb{R}}\int_{\mathbb{R}^d}J(z,r)\mu(0,0)\mu(-z,-r)\nonumber\\
&&\cdot\Big|\nabla u(x,t)\Big(-z+\zeta_{-z,-r}(0,0,\omega)\Big)\Big|^{p-2}\nonumber\\
&&\cdot[-z\zeta_{-z,-r}(0,0,\omega)+z^{2}]\cdot\nabla\nabla u(x,t)dzdr. \label{Equ4.8}
\end{eqnarray}
\end{theorem}

\section{Proof of Lemma \ref{lem3.4}}\label{sec5}
\setcounter{equation}{0}

Let's decompose three parts for the non-local operator $L^\varepsilon$ such that
\begin{eqnarray}\label{Equ5.1}
&&L^\varepsilon v(x,t)=\frac{1}{\varepsilon^{d+4}}\int\int_{\mathbb{R}^{d}\times\mathbb{R}}J(\frac{x-y}{\varepsilon},\frac{t-s}{\varepsilon^2})\nu(\frac{x}{\varepsilon},\frac{t}{\varepsilon^2}, \frac{y}{\varepsilon}, \frac{s}{\varepsilon^2})\nonumber\\
&&\qquad\qquad\cdot\Big(v(y,t)-v(x,t)+v(y,s)-v(y,t)\Big)dyds\nonumber\\
&&=\frac{1}{\varepsilon^{d+2}}\int\int_{\mathbb{R}^{d}\times\mathbb{R}}J(\frac{x-y}{\varepsilon},r)\nu(\frac{x}{\varepsilon},\frac{t}{\varepsilon^2}, \frac{y}{\varepsilon}, \frac{t-\varepsilon^2 r}{\varepsilon^2})\Big(v(y,t)-v(x,t)\Big)dydr\nonumber\\
&&+\frac{1}{\varepsilon^{2}}\int\int_{\mathbb{R}^{d}\times\mathbb{R}}J(z,r)\nu(\frac{x}{\varepsilon},\frac{t}{\varepsilon^2}, \frac{x-\varepsilon z}{\varepsilon}, \frac{t-\varepsilon^2 r}{\varepsilon^2})\Big(v(x-\varepsilon z,t-\varepsilon^2 r)-v(x-\varepsilon z,t)\Big)dzdr\nonumber\\
&&=\frac{1}{\varepsilon^{d+2}}\int\int_{\mathbb{R}^{d}\times\mathbb{R}}J(\frac{x-y}{\varepsilon},r)\nu(\frac{x}{\varepsilon},\frac{t}{\varepsilon^2}, \frac{y}{\varepsilon}, \frac{t}{\varepsilon^2}-r)\Big(v(y,t)-v(x,t)\Big)dydr\nonumber\\
&&-\partial_t v(x,t)\int_{\mathbb{T}^d \times \mathbb{T}}\int_{\mathbb{R}^d\times \mathbb{R}}J(z,r)\nu(\xi,q,\xi-z,q-r)rdzdrd\xi dq +f_\varepsilon(v)=0,
\end{eqnarray}
where use the Taylor expansion and $f_\varepsilon(v)=o(\varepsilon) P(\partial_{t},\partial_{tt},\partial_t \nabla_x\cdots)v\in L^1((0,T)\times L^2(\mathbb{R}^d))$ when $v\in L^2((0,T)\times\mathbb{R}^d)$, $P$ is a combined operator containing higher order space-time derivatives. We introduce a more general Brezis-Pazy Theorem.
\begin{lemma}\label{lem5.1}\cite[Theorem 4.1]{Pazy} Let $A^t_\varepsilon$ be m-accretive in $X, x_\varepsilon \in$ $\overline{D\left(A^t_\varepsilon\right)}, 0\leq\varepsilon<1$ and $f_\varepsilon\in L^1(0,T; X)$ for $n=1,2, \ldots$. Let $u_\varepsilon$ be the mild solution of $\partial_t u_\varepsilon+ A^t_\varepsilon u_\varepsilon \ni f_\varepsilon$ in $[0, T], u_\varepsilon(0)=x_\varepsilon, {D}_\varepsilon =D(A_\varepsilon^0)$.
If $f_\varepsilon \rightarrow f_{0}$ in $L^1(0, T; X)$, $x_\varepsilon \rightarrow x_{0}$ as $\varepsilon \rightarrow 0$ and
$$
\lim_{\varepsilon \rightarrow 0}\left(I+\lambda A^t_\varepsilon\right)^{-1}z=\left(I+\lambda A_{0}\right)^{-1}z, \quad 0\leq t\leq T,
$$
for $z \in Q, 0<\lambda<\lambda_0,$ where $Q=\left(\cap_{1>\varepsilon\geq 0} \bar{D}_\varepsilon\right)\cap\hat{D}, \hat{D}=\hat{D}(A^0)$. Let $U^\varepsilon(t)$ be the evolution operator on $\bar{D}_\varepsilon$ associated with $A_\varepsilon^t$. For each $A_\varepsilon^{t}$  and all $z \in Q$, with $Q$ dense in $X$, then
$$
\lim_{\varepsilon \rightarrow 0} U^\varepsilon(t)z =U(t)z,
$$
that is, $\lim_{\varepsilon\rightarrow 0}u_\varepsilon(t)=u_{0}(t)$ uniformly on $[0,T].$
\end{lemma}
Here we take
\begin{eqnarray}
&&A^t_\varepsilon:=\frac{1}{\varepsilon^{d+2}}\int_{\mathbb{R}^{d}}J_1(\frac{x-y}{\varepsilon})\hat{\nu}(\frac{x}{\varepsilon}, \frac{y}{\varepsilon},\frac{t}{\varepsilon^2})\Big(v(y,t)-v(x,t)\Big)dy,\label{Equ5.2}\\
&&\int_{\mathbb{R}}J(\frac{x-y}{\varepsilon},r)\nu(\frac{x}{\varepsilon},\frac{t}{\varepsilon^2}, \frac{y}{\varepsilon},\frac{t}{\varepsilon^2}-r)dr:=J_1(\frac{x-y}{\varepsilon}) \hat{\nu}(\frac{x}{\varepsilon}, \frac{y}{\varepsilon},\frac{t}{\varepsilon^2}),\label{Equ5.3}
\end{eqnarray}
and we direct give the result in \cite[Remark 3.3]{Piatnitski2022}.
\begin{lemma}\label{lem5.2}\cite{Piatnitski2022} Let $M_l(J):=\displaystyle{\int_{\mathbb{R}^d}|x|^lJ_1(x)dx<\infty}$ for $l=1,2,l_0$ with $l_0$ that satisfies the inequalities $2<l_0<3$. Then we have
$$
\left\|\left(\lambda A_\varepsilon+I\right)^{-1}-\left(\lambda A_0+I\right)^{-1}\right\| \leqslant C \varepsilon^{l_0-2}, \quad \varepsilon>0,
$$
here the constant $C$ depends on $J_1, \mu$ and $l_0$.
\end{lemma}
{\bf The proof of Lemma \ref{lem3.4}} By means of the moments $\displaystyle{\int_{\mathbb{R}}\int_{\mathbb{R}^d}|x|^kJ(x,t)dxdt<+\infty}$ for any $k\geq0$ and (\ref{Equ5.1}), Lemmata \ref{lem5.1}-\ref{lem5.2}, we have $||\varepsilon^\gamma \Sigma^{\varepsilon}||_{L^{2}(\mathbb{R}^{d}\times(0,T))}\leq C \varepsilon^{\beta}, 1 \geq \beta>\gamma>0$. We complete the Lemma \ref{lem3.4}. $\Box$

\section{Proof of exponential convergence in Lemma \ref{lem3.2}}\label{sec6}
\setcounter{equation}{0}

Let $\frac{x}{\varepsilon},\frac{y}{\varepsilon}$ will operate on the torus $\mathbb{T}^d=\Omega$ and $u(x,t)=0$ in $ \Omega_0 { \times \mathbb{ R}^+, \Omega_0 \subset \mathbb{R}^d \setminus \Omega }$. We want to obtain the decay rate in $L^{2}(\Omega)$. Since the  decomposition of $L^\varepsilon$ in (\ref{Equ5.1}), we consider the abstract structure equation
\begin{eqnarray}\label{Equ6.1}
\partial_t u=\int_{\Omega\cup\Omega_0}J(x-y)(u(y,t)-u(x,t))dy+\varepsilon\hat{f}_\varepsilon(u(x,t)).
\end{eqnarray}
Multiplying  the equation by $u$, and integrate to obtain
\begin{eqnarray}\label{Equ6.2}
\frac{1}{2}\partial_t \int_{\Omega} {u^{2}(x, t)} d x &=&\int_{\Omega \cup \Omega_0} \int_{\Omega \cup \Omega_0} J(x-y)(u(y, t)-u(x, t)) u(x, t) d y d x \nonumber\\ &+&\int_{\Omega \cup \Omega_0} \int_{\Omega \cup \Omega_0}\varepsilon \hat{f}_\varepsilon(u(x,t))u(x,t)dx\nonumber\\
&\leq& -\frac{1}{2} \int_{\Omega \cup \Omega_0} \int_{\Omega \cup \Omega_0} J(x-y)|u(y, t)-u(x, t)|^2 d y d x+\varepsilon f_{1\varepsilon}(t).
\end{eqnarray}

On the other hand, using the nonlocal Poincare inequality $\lambda_{1}$ in \cite[Theorem 2.5]{Rossibook} we get
$$
\partial_t \int_{\Omega} \frac{u^{2}(x, t)}{2} d x \leq-C \lambda_{1}\int_{\Omega}|u(x, t)|^{2} d x+C\varepsilon f_{1\varepsilon}(t),
$$
applying the Gronwall inequality
$$
x(t)\leq a(t)+b\int_0^tx(s)ds\Rightarrow x(t)\leq a(t)+b\int^t_0a(s)e^{b(t-s)}ds,
$$
we have exponential decay of the $L^2$ norm
\begin{eqnarray}\label{Equ6.3}
\|u\|^2_{L^{2}(\Omega)}&\leq&\left\|u_0\right\|^2_{L^{2}(\Omega)} e^{-C \lambda_{1}\left(\Omega_0\right) t}+C\varepsilon \int_0^t f_{1\varepsilon} (s)ds -C^2 \varepsilon \lambda_{1} \int_0^t \int_0^s f_{1\varepsilon} (r)dr e^{ -C\lambda_{1}(t-s)}ds\nonumber\\
&\to& \left\|u_0\right\|^2_{L^{2}(\Omega)} e^{-C \lambda_{1}\left(\Omega_0\right) t} \quad \text{as}~ \varepsilon \to 0.
\end{eqnarray}
Equation (\ref{Equ3.1}) has the same exponential decay as nonlocal evolution equation in \cite[Theorem 1.1]{Rossi2015} from above. The inequality $\mathfrak{P}$ and $\chi^T $ converges exponentially to $\chi$ in (\ref{Equ3.12}) are also naturally satisfied.

\section{Construct numeric formats}\label{sec7}
\setcounter{equation}{0}

\subsection{Finite difference discretization}
For $p=2,$ let the positive constant $\nu=1$ and $\delta=\displaystyle{\iint_{\mathbb{R}^{d}\times\mathbb{R}}J(z,r)rdzdr.}$ We take $\mathbf{V}^\varepsilon$ as an approximation solution of (\ref{Equ2.1}) and satisfies that
\begin{eqnarray*}
\left\{\begin{array}{l}
\delta\partial_t\mathbf{V}^\varepsilon(x,t)-\frac{1}{\varepsilon^{d+2}}\displaystyle{\iint_{\mathbb{R}^{d}\times\mathbb{R}}}
J(\frac{x-y}{\varepsilon},r)\Big(\mathbf{V}^\varepsilon(y,t)-\mathbf{V}^\varepsilon(x,t)\Big)dydr=0, \\
\mathbf{V}^\varepsilon(x,0)=f(x), \end{array}\right.
\end{eqnarray*}
where the numerical solution of $\mathbf{V}^\varepsilon$ can be obtain in \cite{Rossinumerical}.

Since the complexity of the calculation and the limitations of storage space, for efficient computation, a sparse matrix representation of the kernel is generated by only evaluating it on nearest neighbors of each of $(x_i,t_h)$. $A_j $ and $A_k=\{i,h$ such that $|x_i-x_j|\leq S, |t_h-t_k|\leq S_1\}$ and $S\times S_1$ is the support of $J(x,t)$, $h_j$ and $\tau_k$ are the maximum unit intervals for adjacent areas of $(x_i,t_h), 1\leq i \leq m, 1\leq h\leq n.$ We consider the space-time discrete numeric format
\begin{eqnarray}\label{Equ7.1}
\left\{\begin{array}{l}
\displaystyle{\hat{\mathbf{L}}^\varepsilon \mathbf{V}^\varepsilon (x_i,t_h)=\frac{1}{\varepsilon^{d+4}} \sum_{j\in B_i}\sum_{k\in B_{h}}J(\frac{x_i-x_j}{\varepsilon},\frac{t_h-t_k}{\varepsilon^2})  \Big(\mathbf{\hat{V}}^\varepsilon(x_j,t_k)-\mathbf{V}^\varepsilon(x_i,t_h)\Big) h_j\tau_k, }\\
\mathbf{\hat{V}}^\varepsilon(x,t)=f(x), t\leq 0, \end{array}\right.
\end{eqnarray}
and compute the error
\begin{eqnarray}\label{Equ7.2}
\mathbf{e}=||\hat{\mathbf{L}}^\varepsilon \mathbf{V}^\varepsilon||_2=\Big(\sum_{i\in A_i}\sum_{h\in A_h}(\hat{\mathbf{L}}^\varepsilon \mathbf{V}^\varepsilon (x_i,t_h))^2 \mathbf{h}_i \hat{\tau}_h\Big)^{\frac{1}{2}},
\end{eqnarray}
where $\mathbf{h}_i$ and $\hat{\tau}_h$ are the sizes of $(x,t)\in D\times (0,T), 2-$norm represents $L^2(\mathbb{R}^d\times (0,T))$, $D$ is a bounded set in $\mathbb{R}^d.$ Actually, we claim that $||u^\varepsilon-u||_2\leq C\varepsilon^\epsilon,$ and $\epsilon$ is arbitrarily close to $1.$

From \cite[Theorem 3.1]{Harlim}, we can deduce that
\begin{eqnarray}\label{Equ7.3}
  ||\mathbf{V}^\varepsilon-u||_2\leq C\varepsilon^\epsilon.
\end{eqnarray}
In fact $\chi_1,\chi_2=0$ when $\nu=1,$ since symmetry of kernel function $J(\cdot,t)$ and (\ref{Equ2.10}), we easily conclude that
\begin{eqnarray}\label{Equ7.4}
  ||\mathbf{V}^\varepsilon-u^\varepsilon||_2\leq C\varepsilon^\epsilon.
\end{eqnarray}

The assertion is true and (\ref{Equ7.1}) is an effective way to test numerical approximation solution $\mathbf{V}^\varepsilon$. Obviously, the above approach will not work under case of $\mu$  is not constant. Therefore, we need to find a way to construct nonlocal numerical correctors in next subsection.

\subsection{Finite element approximations}

First give the notations $\mathbf{V}=H_0^1 (D).$ By $\#$ we denote spaces of periodic functions,
 assuming that
\begin{eqnarray*}
  \nu(x,t,y,s)=\mu(x,t,\frac{x}{\varepsilon},\frac{t}{\varepsilon^2})\mu(x,t,\frac{y}{\varepsilon},\frac{s}{\varepsilon^2})
\end{eqnarray*}
and $p=2$, we generalize numerical approach \cite{Hong} to non-local models.

We denote by $\mathbf{V}_1$ the space $L^2\left(D \times(0,1), \mathbf{V}_{\#}\right)$. Let finite element spaces are  $\mathbf{V}^L \subset \mathbf{V}$ and $\mathbf{V}_1^L \subset L^2(D \times$ $\left.(0,1),\mathbf{V}_{\#}\right) \bigcap L^2\left(D, H_{\#}^1\left((0,1), \mathbf{V}_{\#}^{\prime}\right)\right)$.  $N \in \mathbb{N}$, let $\Delta t=T / N$. We consider the time sequence $0=t_0<t_1<\ldots<t_N=T$ where $t_m=m \Delta t$. Let $g^L \in \mathbf{V}^L$ be an approximation of $g$. Let $t_{m+1 / 2}=t_m+\Delta t / 2$.   We consider the problem with the help of Lemma 2.1: Find $U_{0, m} \in \mathbf{V}^L, U_{1, m} \in \mathbf{V}_1^L$ for $m=1, \ldots, N$ so that
\begin{eqnarray}\label{Equ7.5}
\left\langle B, \psi_0\right\rangle_H+\int_D \int_0^1 \int_Y C \psi_1dy d \tau d x +\int_D \int_0^1 \int_Y A\psi_0 d y d \tau d x=0,
\end{eqnarray}
where
\begin{eqnarray}\label{Equ7.6}
A&=&\int_{\mathbb{R}}\mu(t_{m+1/2},x,y,\tau)\int_{\mathbb{R}^d}J(z,r)\mu(t_{m+1/2},x,y-z,\tau-r){\hspace{-1cm}}\nonumber\\
&&\cdot\Big[-z\otimes \nabla_{x}\frac{U_{1, m+1}+U_{1, m}}{2}(x,t,y-z,\tau-r)+\frac{1}{2}z^{2}\cdot\nabla_{xx} \frac{U_{0, m}+U_{0, m+1}}{2}\Big]dzdr,\nonumber\\
B&=&-\int_{\mathbb{R}}\mu(t_{m+1/2},x,y,\tau)\int_{\mathbb{R}^d}J(z,r)\mu(t_{m+1/2},x,y-z,\tau-r)r\frac{U_{0, m+1}-U_{0, m}}{\Delta t}dzdr,\nonumber\\
C&=&\int_{\mathbb{R}}\mu(t_{m+1/2},x,y,\tau)\int_{\mathbb{R}^d}J(z,r)\mu(t_{m+1/2},x,y-z,\tau-r)\Big[-z\cdot\nabla_{x} \frac{U_{0, m}+U_{0, m+1}}{2}\nonumber\\
&&+\frac{U_{1, m+1}+U_{1, m}}{2}(x,t,y-z,\tau-r)-\frac{U_{1, m+1}+U_{1, m}}{2}(x,t,y,\tau)\Big]dzdr,
\end{eqnarray}
for all $\psi_0 \in \mathbf{V}^L$ and $\psi_1 \in \mathbf{V}_1^L$. Let $u_{0, m}=u_0\left(t_m\right)$ and $u_{1, m}=u_1\left(t_m\right)$ correspond to $U_{0,m}$ and $U_{1,m}$. Problem (\ref{Equ7.5}) has a unique solution and further we can also obtain numerical corrector.

\section*{Acknowledgments} The research is supported by the National Natural Science Foundation of China (No. 12171442). The research of J. Chen is partially supported by the CSC under grant No. 202206160033.

\section*{Data availability statement} Data sharing is not applicable to this article as no new data were created or analyzed in this study.

\section*{Declaration of Interest} The authors declare that they have no known competing financial interests or personal relationships that could have appeared to influence the work reported in this paper.

\section*{Authors' Contributions} Junlong Chen carried out the homogenization theory, and Yanbin Tang carried out the reaction diffusion equations and the perturbation theory of partial differential equations. All authors carried out the proofs and conceived the study. All authors read and approved the final manuscript.

\end{document}